\documentclass[preprint,12pt]{elsarticle}

\usepackage{amssymb}
\usepackage{amsthm}
\usepackage{amsmath}
\usepackage{mathtools}
\usepackage{mathrsfs}
\usepackage{algorithm}
\usepackage[algo2e]{algorithm2e}
\usepackage{algpseudocode}
\usepackage{array}
\usepackage{multirow}
\usepackage{listings}
\usepackage{tabu}
\usepackage{booktabs}
\usepackage{enumerate}
\usepackage{lineno}
\usepackage{fullpage}
\usepackage{xcolor}
\usepackage[colorinlistoftodos]{todonotes}
\usepackage{bm}
\usepackage{bbm}
\usepackage[colorlinks=true]{hyperref}
\usepackage{url}
\usepackage{textcomp}
\usepackage{gensymb}
\usepackage{soul}
\usepackage{graphicx}
\usepackage{subfig}
\usepackage{float}
\usepackage{caption}
\usepackage{tikz}

\usepackage{titlesec} 
\titleformat*{\subsection}{\normalsize\bfseries}
\titleformat*{\subsubsection}{\small\bfseries}

\usepackage[version=4]{mhchem}

\hypersetup{colorlinks=true, citecolor=blue, linkcolor=blue, urlcolor=blue}
\makeatletter

\newcommand{\Rmnum}[1]{\expandafter\@slowromancap\romannumeral #1@}
\makeatother

\usepackage{appendix} 
\usepackage{minitoc} 

\newcommand{\bctnet}{BCTNet}
\newcommand{\gbctnet}{GBCTNet}

\journal{Journal of Computational Physics}

\begin{document}

\begin{frontmatter}



\title{An output scaling layer boosts deep neural networks for multiscale ODE systems}


\author[1]{Yuxiao Yi} 
\author[2]{ Weizong Wang }
\author[2,3]{Tianhan Zhang\corref{cor1}} 
\ead{thzhang@buaa.edu.cn}
\author[1,4]{Zhi-Qin John Xu\corref{cor1}} 
\ead{xuzhiqin@sjtu.edu.cn}

\affiliation[1]{organization={Institute of Natural Sciences, School of Mathematical Sciences,  MOE-LSC, Shanghai Jiao Tong University},city={Shanghai},postcode={200240},country={China}}



\affiliation[2]{organization={School of Astronautics, Beihang University},city={Beijing},postcode={100191},country={China}}
\affiliation[3]{organization={Key Laboratory of Spacecraft Design Optimization and Dynamic Simulation Technologies, Ministry of Education}, 
city={Beijing}, 
postcode={102206}, 
country={China}}
\affiliation[4]{organization={Key Laboratory of Marine Intelligent Equipment and System, Ministry of Education},postcode={201100},country={China}} 
\cortext[cor1]{Corresponding authors.}

\begin{abstract}

Simulating complex diffusion–reaction systems is often prohibitively expensive due to the high dimensionality and stiffness of the underlying ODEs, where state variables may span tens of orders of magnitude. 
Deep learning has recently emerged as a powerful tool in scientific computing, achieving remarkable progress in modeling and sampling stiff systems. 
However, data scaling techniques remain largely underexplored, despite their crucial role in addressing the frequency bias of deep neural networks when handling multi-magnitude or high-frequency data.
In this work, we propose the Generalized Box–Cox Transformation (GBCT), a novel nonlinear scaling method designed to mitigate multiscale challenges by rescaling inherent multi-magnitude components toward a more consistent order of magnitude. 
We integrate GBCT into our previous data-driven framework and evaluate its performance against the original baseline surrogate model across six representative scenarios: a 21-species chemical reaction kinetics, a 13-isotope nuclear reaction model, the well-known Robertson problem coupled with diffusion, and practically relevant simulations of two-dimensional turbulent reaction-diffusion systems as well as one- and two-dimensional nuclear reactive flows. 
Numerical experiments demonstrate that GBCT reduces prediction errors by up to two orders of magnitude compared with the baseline model—particularly in the long-term evolution of dynamical systems—and achieves comparable performance with only about one-sixth of the training epochs.
Frequency analysis further reveals that GBCT rescales high-frequency components of the objective function toward lower frequencies to align with the neural network’s natural low frequency bias, thereby boosting training and generalization. The source code to reproduce the results in this paper is available at \url{https://github.com/Seauagain/GBCT}.

\end{abstract}




\end{frontmatter}

\doparttoc 
\faketableofcontents 
\part{} 
\vspace{-2em} 


\section{Introduction}
\label{sec:intro}

Reacting flows play a crucial role in various natural phenomena and engineering applications, including chemical engineering, aerospace engineering, astrophysical nucleosynthesis, biochemistry, and atmospheric environments. 
However, numerical simulation of reacting flows encounters substantial computational challenges rooted in the strong stiffness of the high-dimensional ordinary differential equations (ODEs) governing reaction dynamics \cite{zingale2024}. 
These challenges arise because numerical schemes typically require small time steps to ensure stability and accuracy, which results in prohibitive computational overhead. Therefore, to accelerate reacting flow simulations, a novel method for handling stiff reaction kinetics is urgently needed.


Deep learning has emerged as a powerful tool for scientific computing, providing new insights into accelerating the computation of stiff ODEs. 
Early studies leveraged deep neural networks (DNNs) to model the chemical reaction source terms \cite{christo1996CF, christo1996SC, blasco1998CF, blasco2000CTM, kempf2005PCI, ihme2009PCI}. 
Recent research interest in reacting flows has shifted towards simulations reflecting realistic industrial scenarios, requiring large-scale dynamical equations to model reactions and incorporating turbulence, thereby creating an urgent demand for enhanced accuracy and generalization capability in modern deep learning methods.
The key challenge lies in the intricate multiscale nature of the reaction dynamics, simultaneously posing difficulties for modeling, training and sampling. 
To address the challenges, many attempts have emerged, focusing on modeling approaches, training algorithms and sampling methods. 
The emerging paradigms of physics-informed neural networks (PINNs) \cite{raissi2019JCP} and neural ordinary differential equations (neural ODEs) \cite{chen2018ANIP} gain popularity in scientific computing, inspiring researchers to model reaction dynamics. 
Owoyele and Pal \cite{owoyele2022EA} proposed a novel framework ChemNODE based on neural ODEs to model hydrogen reaction kinetics. 
To address the training imbalance caused by multiscale nature, Ji et al. \cite{ji2021JPCA} introduced stiffPINN embedding quasi-steady-state assumptions (QSSA) while Weng et al. \cite{weng2022JPCA} designed MPINN utilizing separate neural networks to reduce stiffness. Florio et al. \cite{deflorio2022CIJN} combine PINNs with the theory of functional connections and extreme learning machines to tackle stiff chemical kinetics.
Subsequent works combined with neural operators \cite{lu2021NMI, li2020APA} advanced the field by validating the performance of DeepONet variants \cite{goswami2024CMAM} and the extended FNO \cite{weng2025CF} in more challenging turbulent reacting flows. 
Another significant branch of development focuses on enhancing the generalization 
capability of data-driven models in complex scenarios through sampling. 
Ideal training data is desired to be sufficient and representative, capable of covering samples from unseen test scenarios. 
To this end, Wan et al. \cite{wan2020CF} generated data from low-cost canonical problem simulations and enhanced it with random perturbations. 
Ding et al. \cite{ding2021CF} and Readshaw et al. \cite{readshaw2023CF} developed a hybrid database combining random and simulated data.
Zhang et al. \cite{zhang2022CF} proposed a multi-scale sampling algorithm combined with data scaling to mitigate multi-scale effects in data, achieving high accuracy of surrogate models for 9-dimensional hydrogen kinetics in reactive flow simulations. 
Yao et al. \cite{yao2025CPC} developed a novel sampling method that enhances random data through temporal evolution, significantly improving data quality and model generalization. 
The corresponding algorithm has been extensively validated for robustness and accuracy across different dynamic systems \cite{zhang2025AJ}. 
To address the lack of physical conservation laws in data-driven models, Wang et al. \cite{wang2025CF} strictly enforced element, mass, and energy conservation laws into the prediction of reaction source terms based on the previous sampling methods \cite{zhang2022CF, yao2025CPC}. These works bridge the gap between data distribution, physical constraints and model performance, providing new insights for the field.

Prior research has largely concentrated on modeling and data sampling while the area of data scaling remains relatively unexplored, which holds potential for addressing challenges in modeling.
Deep neural networks exhibit an inherent tendency known as the frequency principle \cite{Xu2020CiCP, xu2025CAMC, xu2019NIP} or spectral bias \cite{rahaman19a} during training, where they prioritize capturing the low-frequency components of the target function. 
However, the multiscale nature of these systems entails a broad frequency spectrum, complicating training and hindering convergence. Data scaling techniques are designed to transform high-frequency components into low-frequency ones, thereby enabling models to more effectively discern the intrinsic patterns of dynamical systems and enhance generalization. Kim and Ji et al. \cite{kim2021CIJN} proposed feature scaling by linearly adjusting the inputs, labels, and loss functions of neural ODEs, ensuring that the contributions of various species to the total loss function were normalized to the same magnitude. The efficacy of these scaling techniques was demonstrated using the classic stiff Robertson problem. Zhang et al. \cite{tz2021} initially employed the nonlinear Box-Cox transformation (BCT) to deal with the chemical species concentrations as inputs of DNNs. 
Extensive subsequent studies employing BCT have further confirmed its effectiveness in addressing multiscale challenges \cite{weng2025CF,yao2025CPC, zhang2025AJ, wang2025CF,lian2023AA,saito2023AECS,lei2024ANE,wu2024EA,zhang2024PF, li2025,rao2025CF}. 
The prevalence of BCT in reaction dynamics modeling has attracted our interest in exploring data transformation.
However, BCT’s non-negativity constraint limits its ability to handle sign-changing multiscale data, which often serve as the outputs of neural networks. 
Therefore, a new data scaling method is needed.

In this work, we propose a novel scaling method termed the Generalized Box-Cox Transformation (GBCT). 
This method effectively mitigates the multiscale effects present in sign-changing training data of DNNs. 
We evaluate this method across a diverse set of benchmarks, including chemical and nuclear reaction kinetics ODEs, the Robertson diffusion problem, as well as challenging chemical and nuclear reactive flow simulations.
Integrating GBCT into our previously developed data-driven framework \cite{yao2025CPC}, we assess its impact on model accuracy, long-term evolution stability, and error accumulation in both pure dynamics and complex turbulent diffusion-reaction settings. 
Furthermore, we experimentally analyze the influence of GBCT on training efficiency, revealing that it accelerates the neural network training process. 

The remainder of the paper is organized as follows. 
Sec.~\ref{sec:problem} introduces a representative multiscale problem, reaction kinetics, which forms the central focus of this study.
Sec.~\ref{sec:method} outlines the methodology on surrogate modeling integrated with GBCT. 
Sec.~\ref{sec:gbct} presents the underlying inspiration of GBCT.
Sec.~\ref{sec:result} demonstrates extensive benchmark results, highlighting the improvements brought by GBCT in terms of model accuracy, versatility, and long-term stability.
Sec.~\ref{sec:analysis} further examines the impact of GBCT on training efficiency and generalization capability, and analyzes its effectiveness through a frequency perspective.
Finally, Sec.~\ref{sec:conclusion} concludes the paper.

\section{Problem Description}
\label{sec:problem}

This work focuses on a widely applied class of reaction kinetics models based on molecular collision.
A general elementary reaction involving $N_s$ species can be represented as:
\begin{align}
    \sum_{k=1}^{N_s} \alpha_k S_k \rightarrow \sum_{k=1}^{N_s} \beta_k S_k,
\end{align}
where $S_k$ is the $k$-th species, $\alpha_k$ and $\beta_k$ are the stoichiometric coefficients of reactants and products, respectively. These coefficients are integers or fractions that balance the reaction, ensuring the conservation of elements across the reaction. The molar production rate of the $i$-th species, $\dot{\omega}_i$ (moles per unit volume per second), is defined as: 
\begin{align}
    \frac{\mathrm{d} [X_i] }{\mathrm{d} t} = \dot{\omega}_i = (\beta_i-\alpha_i)\cdot k_f \cdot \prod_{k=1}^{N_s} [X_k]^{\alpha_k},
    \label{eq:general_rate}
\end{align}
where $[X_i]$ denotes the molar concentration (moles per unit volume) of species $S_i$, $k_f$ is the forward reaction rate constant and $(\beta_i-\alpha_i)$ is the net stoichiometric coefficient. The exponents $\alpha_k$ introduce a nonlinear term, amplifying the sensitivity of the reaction rate to changes in species concentrations. Note that the molar concentration and mass fraction are related by $ [X_i] = \tfrac{\rho Y_i}{M_i}$, where $M_i$ is the molecular weight and $\rho$ is the density.
Hence, the species production rate can also be equivalently expressed in the form based on the mass fraction:
\begin{equation}
    \frac{\mathrm{d} Y_i }{\mathrm{d} t } = (\beta_i-\alpha_i)\frac{M_i}{\rho} \cdot k_f\cdot  \prod_{k=1}^{N_s} \left(\frac{\rho Y_k}{M_k} \right)^{\alpha_k}.
    \label{eq:ode_rate_mass}
\end{equation}
The unified formula exhibits remarkable generality, extending beyond fuel chemical
kinetics to diverse physical and engineering processes. The differences among reaction systems mainly lie in the specific form of reaction rate constant $k_f$. Building upon the general form, we will briefly introduce several representative reaction kinetics from distinct application fields.

\subsection{Simplified Reaction Kinetics With Constant $k_f$}

The Robertson problem (BOBER) \cite{robert1966} is a classical stiff ODE systems consisting of three species and three reactions with widely separated rate constants. We consider an elementary reaction in the Robertson problem, $ B + B \rightarrow B + C $, where the mass fraction of species $B$ is denoted by $y_2$. The governing equation for its temporal evolution is  
\begin{align}
\frac{\mathrm{d} y_2}{\mathrm{d} t} = -k_2 y_2^2.
\end{align}
Mathematically, we take $\rho=M_i=1$ for simplicity which makes this equation consistent with the general form in Eq.~(\ref{eq:ode_rate_mass}). In the ROBER problem, the rate constants are constant, and the rates of the three elementary reactions span nine orders of magnitude, ranging from $\mathcal{O}(10^{-2})$ to $\mathcal{O}(10^7)$, resulting in strong stiffness within the ODE system.


\subsection{Nuclear Reaction Kinetics}

In nuclear reactions, the evolution of nucleus A involved in a typical two-body reaction, $A+B\rightarrow C+D$, takes the form \cite{smith2023A}:
\begin{align}
  \frac{\mathrm{d} Y_1}{\mathrm{d} t}= -1 \cdot \frac{M_1}{\rho} \cdot  N_{\mathrm{A}}\langle\sigma v\rangle \cdot \left( \frac{\rho Y_1}{M_1}\right) \left( \frac{\rho Y_2}{M_2}\right),
\end{align}
where $Y_1$ and $Y_2$ denote the mass fraction, $M_1$ and $M_2$ denote the nucleus weight of the species $A$ and $B$, and $N_{\mathrm{A}}$ is Avogadro's constant. The rate constant is empirically modeled as a function of temperature $T$:
\begin{align}
     k_f := N_{A}\langle\sigma v\rangle = \exp\bigg(a_{0}+\sum_{i=1}^{5}a_{i}T_{9}^{(2i-5)/3}+a_{6}\log T_{9}\bigg),
\end{align}
with $T_9 = T / (10^9 K)$ and $a_i$ as fitted coefficients. The evolution of the governing ODE system of nuclear reaction depends on the temperature $T$, density $\rho$ and the mass fractions $Y_i$. The temperature $T$ can be calculated through the energy conservation and $\rho$ is updated by the equation of state (EOS).

\subsection{Elementary Chemical Reaction Kinetics}

The third representative reaction system is detailed reaction model consisting of elementary chemical reactions. The forward reaction rate constant $k_f$ is specific to each reaction and primarily depends on temperature. In some cases, it also varies with pressure, a phenomenon known as a pressure-dependent reaction. The form of the reaction rate $k_f$ is given by the Arrhenius law 
\begin{equation}
    k_f=Ae^{-E_a/{RT}},
\end{equation}
where $A$  is the pre-exponential factor, $E_a$ is the activation energy, $R$ is the universal gas constant, and T is the absolute temperature. Mathematically, this exponential form captures the rapid increase in reaction rate with temperature, a critical feature in stiff systems. Similar to nuclear reaction kinetics, the physical state vector of the ODE system consists of temperature, density and mass fractions.

Although these reaction models mentioned above take different forms, they share similar stiffness or multi-scale feature, caused by both mass fractions and reaction rate constants.
This resulting stiffness, characterized by eigenvalues of the Jacobian matrix spanning multiple orders of magnitude, demands specialized numerical solvers to ensure stability and accuracy. 
Such ODE systems are crucial in computational physics, as they underpin simulations of combustion, nuclear fusion, and material evolution, enabling the design of efficient energy systems, predicting stellar nucleosynthesis, and optimizing biochemical reactors, among other applications.

\subsection{Surrogate Modeling}
\label{sec:surrogate_model}

In conventional reaction-diffusion simulations, the reaction dynamics can be decoupled from transport process through operator splitting, such as Strang splitting \cite{Strang1968}, allowing independent calculations at each grid cell. 
Therefore, a common approach is to model the most time-consuming ODE integration to overcome the computational bottleneck of the simulation. The general form of ODE can be expressed as 
\begin{equation}
\label{eq:generic_ode}
    \frac{ \mathrm{d} \pmb{x}}{\mathrm{d} t } = \pmb{\dot{\omega}} (\pmb{x}, t), 
\end{equation}
where $\pmb{x}:=\pmb{x}(t)\in \mathbb{R}^{d}$ is the state vector usually consisting of species concentrations, density and temperature, $d$ is the dimensionality associated with the number of species, and $\pmb{\dot{\omega}}$ denotes the reaction source term. The temporal evolution of the reaction dynamical system is defined as
\begin{align}
\label{eq:evo_ode}
    \pmb{x}(t+\Delta t) := \mathcal{F}(\pmb{x}(t)) = \pmb{x}(t)+\int_{t}^{t+\Delta t} \pmb{\dot{\omega}}\mathrm{d} t ,
\end{align}
where $\mathcal{F}:\mathbb{R}^{d}\rightarrow \mathbb{R}^{d}$
denotes the strongly non-linear propagation operator advancing the state $\pmb{x}(t)$ to $\pmb{x}(t+\Delta t)$ with a desired time step $\Delta t$.
The choice of a desired time step $\Delta t$ is critical. Stiff ODEs require an extremely small time step in traditional solvers to maintain stability, whereas a larger, user-defined time step $\Delta t$ facilitates efficient long-term predictions, significantly accelerating simulations and alleviating the computational bottleneck posed by the reaction source term.

A prevalent approach to build a surrogate model is leveraging the deep neural network (DNN)  $f_{\pmb{\theta}}$ with learnable parameters $\pmb{\theta}$ to replace the integral. 
The DNN takes the current physical state $\pmb{x}(t)$ as input. Alternatively, instead of predicting the state $\pmb{x}(t+\Delta t)$, the label of DNN can be set as the state change rate, defined as $\pmb{u}(t;\Delta t):= \left(\pmb{x}(t+\Delta t) - \pmb{x}(t) \right)/\Delta t$. The optimal parameters $\pmb{\theta}^{*}$ minimize the loss between the predicted and actual state change rates:
\begin{align}
\pmb{\theta}^{*} = \arg\min\limits_{\pmb{\theta}}~ \mathcal{L} \left( f_{\pmb{\theta}}(\pmb{x}(t)), ~\pmb{u}(t;\Delta t)  \right),
\end{align}
where $\mathcal{L}$ is the loss function, guiding the training process to achieve accurate approximations.

\section{Methodology}
\label{sec:method}

Our study focuses on effectively modeling reaction kinetics which is a representative multiscale problem.
For generic chemical and nuclear reactions , the state vector of the dynamical system can be written as 
\begin{equation}
\label{eq:x}
    \pmb{x}(t): = \left[ T(t),\, \rho(t),\, \pmb{Y}(t) \right], 
\end{equation}
where $\pmb{Y}(t):=\left[ Y_1(t), Y_2(t), \cdots, Y_{N_s}(t) \right]$. 
The components $T(t)$, $\rho(t)$, $Y_i(t)$ and $N_s$ denote temperature, density, mass fraction of the $i$-th species, and the number of reacting species, respectively. 
For other types of reaction systems, the state vector can be modified according to the variables on which the reaction rate constant depends. 
For instance, in the Robertson problem, the state vector depends only on the mass fractions of the species. 
Moreover, $\rho(t)$ can be replaced by pressure $p(t)$, as the two are interchangeable through the equation of state (EOS).
Without loss of generality, the state vector is adopted throughout this paper as a representative example.
 With a slight abuse of notation, we modify the definition of the label $\pmb{u}(t; \Delta t )$ as   
\begin{equation}
\label{eq:u}
    \pmb{u}(t;\Delta t): = \frac{\pmb{Y}(t+\Delta t) - \pmb{Y}(t) }{\Delta t}. 
\end{equation}
In our framework, the time step $\Delta t$ is fixed, and thus the $\pmb{u}(t;\Delta t)$ will be simplified to $\pmb{u}(t)$. The temperature and density at time $t+\Delta t$ are not directly predicted by the DNN. 
Instead, to satisfy the physical constraints, $\tilde{T}(t+\Delta t )$ is calculated through energy conservation equation while $\tilde{\rho}(t+\Delta t )$ is derived from the equation of state (EOS) using the current state $\pmb{x}(t)$ and the DNN-predicted $\tilde{\pmb{Y}}(t+\Delta t)$.
Therefore, accurate prediction of mass fractions is of crucial importance. 

\subsection{Data Transformation}
\label{sec:transform}

During the reaction process, the species mass fractions exhibit typical multi-scale characteristics. 
Across different species or evolutionary time instants, the characteristic scales span multiple orders of magnitude, ranging from $\mathcal{O}(10^{-25})$ to $\mathcal{O}(10^{-1})$.
Due to the  frequency principle of generic neural networks, DNNs struggle to fit $\pmb{u}(t)$ when directly taking the multi-scale terms $\pmb{Y}(t)$ in as inputs.
Therefore, the raw data needs to be appropriately scaled to near $\mathcal{O}(1)$ which ensures consistent magnitudes across distinct species, thereby preventing gradient imbalance during DNN training \cite{wang2021}.
Previous studies have addressed the multi-scale issue in the inputs, but we observed that strong scale effects still exist in the outputs which are often overlooked. 
In the following, we will separately discuss the individual handling of inputs and outputs.

\begin{figure}[htb]
	\centering
\includegraphics[width=0.95\linewidth]{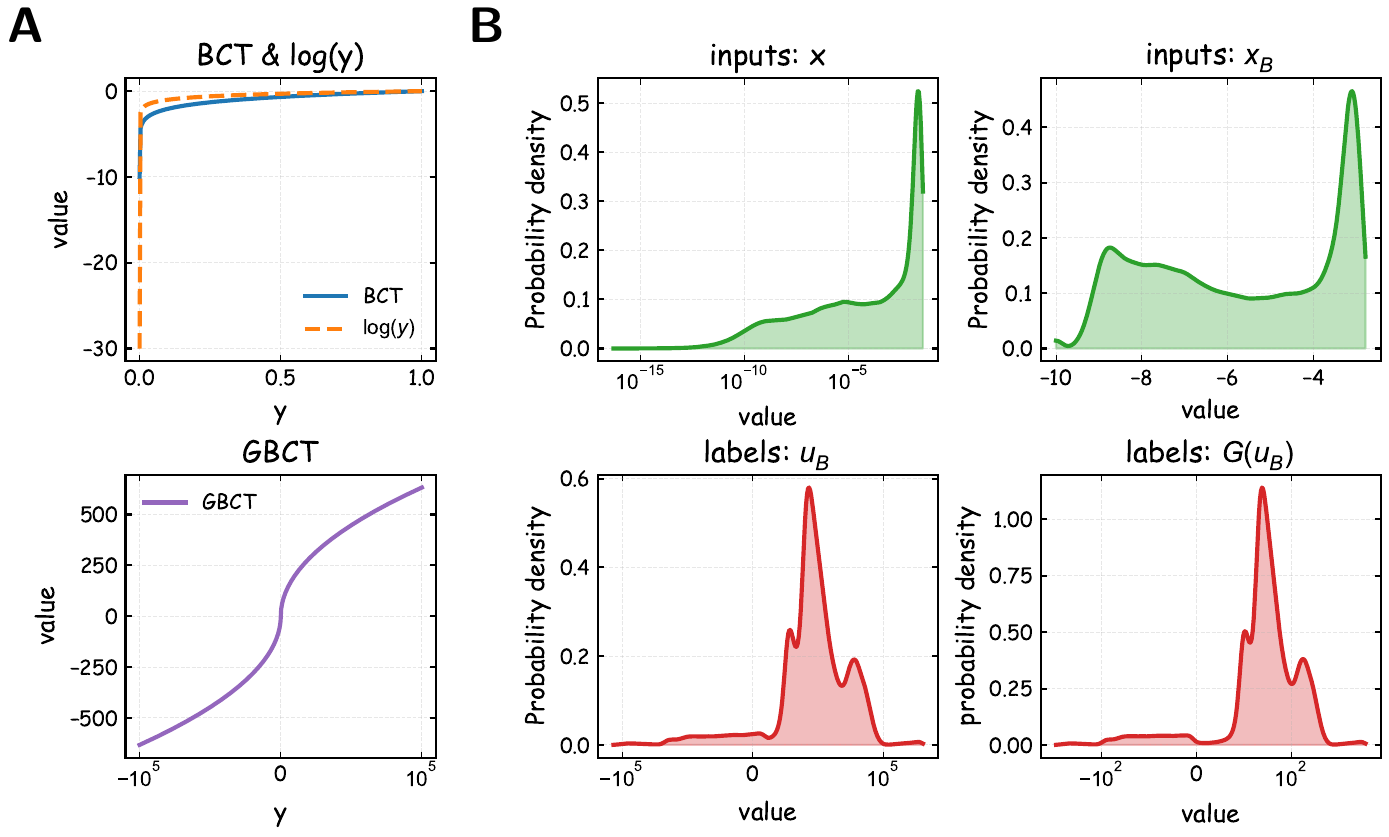}
	\caption{(A) Schematic diagrams of $\log(x)$, BCT and GBCT transformations. (B) Distribution of 
 species \ce{H} mass fraction in inputs and outputs before and after data transformation. The BCT maps $\pmb{x}(t)$ to $\pmb{x}_B(t)$ while the GBCT maps $\pmb{u}_B(t)$ to $G(\pmb{u}_B(t))$ with parameters $\lambda_a = 0.1$ and $\lambda_b = 0.5$.}
	\label{fig:bct_gbct_diagram}
\end{figure}

Considering $Y_i(t) \in [0,1]$ within the inputs, logarithm transformation seems natural but fails to handle scenarios where the mass fractions approach 0, indicating that the logarithmic transformation is ill-suited in our practice. 
To address this issue, the Box-Cox transformation (BCT) is introduced. 
The BCT is a statistical method originally proposed by Box and Cox \cite{Box1964} to stabilize variance of data, and the function is defined as 
\begin{align}
    B(x) = 
\begin{cases}
\displaystyle \frac{x^{\lambda_a} - 1}{\lambda_a}~ & \text{if } \lambda_a \neq 0, \\
\ln(x)~ & \text{if } \lambda_a = 0, 
\end{cases}
\end{align}
where $x$ is the input variable and $\lambda_a$ is a hyper-parameter. 
It is first employed to deal with the concentration of chemical species in the previous literature \cite{tz2021}. 
BCT monotonically maps the non-negative mass fractions within $[0,1]$ to $[-1/\lambda_a, 0]$ when $\lambda_a > 0$. 
Additionally, as shown in the first row of Fig.~\ref{fig:bct_gbct_diagram}(A), BCT converts the low-magnitude concentrations $10^{-p}~(p\geq1)$  to approximately $\mathcal{O}(1)$, thus avoiding the singularity when approaching zero compared with the logarithm function. 
We apply BCT as a data transformation method for mass fractions appearing in Eq.~(\ref{eq:x}) and Eq.~(\ref{eq:u}), and the scaled input and label for DNNs turn out to be 
\begin{align}
\pmb{x}_{B}(t): &= \left[ T(t), \rho(t), B(\pmb{Y}(t)) \right],  \label{eq:xB} \\
\pmb{u}_{B}(t): &= \frac{B(\pmb{Y}(t+\Delta t) ) - B (\pmb{Y}(t) )}{\Delta t} . \label{eq:uB}
\end{align}
The top row of Fig.~\ref{fig:bct_gbct_diagram}(B) displays the distributions of $\pmb{x}(t)$ and $\pmb{x}_B(t)$, with data selected from representative chemical reaction processes for illustration. Obviously, BCT successfully maps species mass fractions spanning more than 10 orders of magnitude (e.g. $10^{-15}\sim10^{-4}$) to near $\mathcal{O}(1)$, significantly mitigating the multi-scale effects. 

Although BCT effectively addresses the multiscale components in the raw input state $\pmb{x}(t)$, as demonstrated in the bottom row of Fig.~\ref{fig:bct_gbct_diagram}(B), the corresponding scaled label $\pmb{u}_B(t)$ still exhibits multiscale nature, ranging from $-10^{5}$ to $10^{5}$. 
Generally, this term $\pmb{u}_B(t)$ tends to $B'(\pmb{Y})\cdot\mathrm{d} \pmb{Y}/\mathrm{d}t$ when $\Delta t$ approaches 0, exhibiting intrinsic multiscale effect essentially introduced by the term $\mathrm{d} \pmb{Y}/\mathrm{d}t$, thereby necessitating additional data scaling method. 
Before describing the details of our method, it is worth clarifying three essential requirements that an ideal transformation is expected to satisfy.
1) The transformation must be monotonic to ensure invertibility during the model inference stage;
2) The transformation is capable of handling negative data with multiscale feature which the BCT fails to address;
3) The transformation should be in simple form, computationally efficient, and should not compromise the model's overall efficiency.

We propose the Generalized Box-Cox Transformation (GBCT) to further process inherent multiscale terms in $\pmb{u}_{B}(t)$, which are persistently overlooked in prior works \cite{zhang2022CF, yao2025CPC}. 
Our method extends the conventional BCT which is limited to non-negative values.
To achieve this, we construct a continuous monotonic function on $\mathbb{R}$ via an odd extension based on BCT.  
Specifically, we first remove the constant term $-1/\lambda_a$ from the original BCT, retaining only the power function. This simplified function is then mirrored about the origin to handle negative inputs. 
The resulting GBCT is defined in the following piecewise form:
\begin{align}
    G(x) = 
    \begin{cases}
        x^{\lambda_b} / \lambda_b, & \text{if } x \geq 0, \\
        -(-x)^{\lambda_b} / \lambda_b, & \text{if } x < 0,
    \end{cases}
    \label{eq:gbct_def}
\end{align}
where $ x $ is the input variable and $ \lambda_b$ is a positive hyper-parameter. 

The term $ \pmb{u}_{B}(t) $ exhibits variations spanning multiple orders of magnitude, encompassing both small-scale fluctuations near the reaction equilibrium and large-scale states during active reaction stages. 
When applied with GBCT, data within the interval $ [-10^0, 10^0] $ are expanded to the order of $ \mathcal{O}(1) $, analogous to the standard BCT. 
In contrast, values outside this range are substantially compressed in magnitude. 
As shown in the bottom row of Fig.~\ref{fig:bct_gbct_diagram}(B), the GBCT effectively rescales $ \pmb{u}_{B}(t) $ from a range spanning over ten orders of magnitude to approximately $[-10^2, 10^2]$.
This scaling operation contributes to a relatively smooth distribution of $ G(\pmb{u}_{B}(t)) $, which benefits the training of deep neural networks by alleviating the difficulties associated with multi-magnitude data in the dynamical systems. 
The inspiration and theoretical analysis behind this idea are further detailed in Sec.~\ref{sec:gbct}.


To assess the effectiveness of GBCT as a simple yet powerful scaling method for reaction dynamics data, we incorporated it into the previously proposed dynamical system modeling framework DeePODE \cite{yao2025CPC} for comprehensive evaluation. 
Since this study centers on the influence of data-level processing on model performance, we preserved the original sampling method, network architecture, and training strategy of DeePODE.
Two distinct neural network models were designed: the baseline model \textbf{\bctnet} and its counterpart \textbf{\gbctnet}.
Both networks receive $\pmb{x}_{B}(t)$ as input.
The unique difference lies in the prediction targets: the baseline model \bctnet{} is trained to predict $\pmb{u}_{B}(t)$, whereas \gbctnet{} utilizes $G\left( \pmb{u}_{B}(t) \right)$ as its label. 
The implementation of GBCT is simple and effective, thereby introducing minimal computational overhead. 
Furthermore, GBCT is agnostic to modeling paradigms, network architectures, and loss function designs, making it a versatile, plug-and-play adaptor module that can be potentially integrated into mainstream architectures such as PINN or DeepONet.

\subsection{Neural Network Structure and Training}

In previous study we developed a multi-scale sampling strategy that randomly samples across the entire reaction state space and employs an evolutionary data augmentation \cite{yao2025CPC}, resulting in high-quality datasets. 
The effectiveness of this sampling method has been extensively validated \cite{wan2020CF,yao2025CPC,zhang2025AJ}. 
We continue to use this reliable strategy to ensure that data quality does not become a bottleneck for network performance in this paper. 
For a fair comparison, \bctnet{} and \gbctnet{} are trained and evaluated on exactly the same dataset. To perform training, the dataset is partitioned into a training set comprising 90\% of the data and a validation set containing the remaining 10\%. The loss function is absolute mean error, i.e., $L_1$ loss. For \bctnet{}, the loss is defined as 
\begin{align}
    L_{\text{BCT}} = \frac{1}{N}\sum_{i=1}^{N}\left\Vert f_{ \pmb{\theta} }( \pmb{x}^{(i)}_B(t) ) -  \pmb{u}^{(i)}_{B}(t)  \right\Vert_{L_1}\ ,
\end{align}
where $\{ (\pmb{x}^{(i)}_B(t), \pmb{u}^{(i)}_{B}(t)  ) \}_{i=1}^{N}$ are the training samples, $N$ denotes the total training size and $f_{ \pmb{\theta} }$ is the DNN model with parameters $\pmb{\theta}$.  In contrast, the loss of \gbctnet{} is defined as 
\begin{align}
    L_{\text{GBCT}} = \frac{1}{N}\sum_{i=1}^{N}\left\Vert f_{ \pmb{\theta} }( \pmb{x}^{(i)}_B(t) ) -  G(\pmb{u}^{(i)}_{B}(t))  \right\Vert_{L_1}.
\end{align}
Note that $\pmb{x}^{(i)}_B(t)$, $\pmb{u}^{(i)}_{B}(t) $ and the function $G(\cdot)$ are defined in Sec.~\ref{sec:transform}. 
We employ two identical fully-connected neural networks, each consisting of three hidden layers with 1600, 800, and 400 neurons, respectively, and use the GELU activation function. 
The Adam optimizer is adopted for training. We implement Z-score normalization (i.e., subtract the mean and divide by the standard deviation) across all input and label data to enhance the convergence speed of DNNs.
A two-stage training strategy is utilized: the first stage runs for 2500 epochs with a batch size of 1024 and a learning rate of $10^{-4}$; the second stage also runs for 2500 epochs, in which the batch size is increased by a factor of $2^8$ and the learning rate is reduced to $10^{-5}$.  
This approach proves effective in enhancing both the accuracy and generalization ability of the neural networks.

\section{Generalized Box-Cox Transformation}
\label{sec:gbct}

In this section, we further provide a step-by-step theoretical analysis to elucidate the intuition and inspiration behind our proposed GBCT method. We aim to answer three key questions: (1) How do we identify that label terms still contain multi-scale components that are often overlooked? (2) Why is GBCT necessary? and (3) What specific problem does GBCT solve?

Practical reaction systems consist of numerous elementary reactions and species, resulting in highly nonlinear terms which are not suitable for theoretical analysis. To clearly describe our idea, we begin our analysis with an irreversible elementary reaction involving $N_s$ species. 
As mentioned before, the original inputs and labels of the DNNs are denoted as $\pmb{x}(t)$ and $\pmb{u}(t)$, respectively.
When $\Delta t \to 0$, $\pmb{u}(t)$ approximates $\mathrm{d} \pmb{Y}(t)/\mathrm{d} t$. We therefore consider the limit form 
$\mathrm{d} \pmb{Y}(t)/\mathrm{d} t$ instead of $\pmb{u}(t)$, since the former is analytically tractable. 
The corresponding inputs and labels at this limit are expressed as follows:
\begin{align}
\pmb{x}(t) &=
\begin{bmatrix}
T(t) \\
\rho(t) \\
Y_1(t) \\
\vdots \\
Y_{N_s}(t)
\end{bmatrix},\quad
\pmb{u}(t) =
\begin{bmatrix}
\displaystyle \frac{\mathrm{d} Y_1}{\mathrm{d} t} \\
\vdots \\
\displaystyle \frac{\mathrm{d} Y_{N_s}}{\mathrm{d} t}
\end{bmatrix}.
\end{align}
As described in Sec.~\ref{sec:problem}, the reaction rate of the $i$-th species in an elementary reaction is generally expressed as:
\begin{align}
    \frac{\mathrm{d} Y_i }{\mathrm{d} t } &= (\beta_i-\alpha_i)\frac{M_i}{\rho} \cdot k_f\cdot  \prod_{k=1}^{N_s} \left(\frac{\rho Y_k}{M_k} \right)^{\alpha_k}, \\
    & \overset{\triangle}{=} C \cdot k_f\cdot \prod_{k=1}^{N_s}[Y_k]^{\alpha_k},
    \label{eq:general_rate_mass}
\end{align}
where $C$ is a time-independent coefficient absorbing the stoichiometric coefficients $\alpha_k$, density $\rho$, and molecular weight $M_i$. Due to the multiscale nature of the mass fraction $Y_k$ in the inputs, we apply BCT to address it. 
After scaling, the transformed inputs and labels are given by:
\begin{align}
\pmb{x}_B(t) &=
\begin{bmatrix}
T(t) \\
\rho(t) \\
B(Y_1) \\
\vdots \\
B(Y_{N_s})
\end{bmatrix},\quad
\pmb{u}_B(t) =
\begin{bmatrix}
\displaystyle \frac{\mathrm{d} B(Y_1) }{ \mathrm{d} t } \\
\vdots \\
\displaystyle \frac{\mathrm{d} B(Y_{N_s}) }{ \mathrm{d} t }
\end{bmatrix}.
\end{align}
We then further expand the term $\mathrm{d} B(Y_{i}) / \mathrm{d} t $ to analyze the underlying multi-scale behavior, which can be denoted as:
\begin{equation}
\label{eq:bct_ode}
   \frac{\mathrm{d} B(Y_i) }{\mathrm{d} t } =C \cdot k_f\cdot B'(Y_i)\cdot \prod_{k=1}^{N_s}[Y_k]^{\alpha_k},
\end{equation}
where $B'(Y_i)$ denotes the first order derivative of $B(x)$ with respect to $Y_i(t)$. 
Although the transformed $B(Y_i)$ eliminates multiscale features in the input, the term $\pmb{u}_B(t)$ still intrinsically contains multiscale components $k_f$ and $\prod_{k=1}^{N_s}[Y_k]^{\alpha_k}$. 
This key observation motivates the necessity of an additional transformation. 

An ideal transformation applied to $\mathrm{d}B(Y_i)/\mathrm{d}t$ should rescale both the high-magnitude rate constant $k_f$ and the low-magnitude $Y_k$ to values near the order of $\mathcal{O}(1)$, while preserving monotonicity, accommodating negative inputs, and remaining easy to implement.
Intuitively, the logarithmic function seems to be a reasonable choice.
However, it suffers from singularities when $Y_k$ approaches zero. 
On the other hand, the BCT, as a variant of the power function, cannot accept negative inputs.
For instance, when $\lambda_a = 0.5$, BCT is equivalent to taking a square root. Therefore, both the logarithmic function and the BCT are not suitable for the transformation.

To overcome the limitations of the logarithmic function and the BCT, we propose Generalized Box-Cox Transformation (GBCT) to achieve data pre-processing for the label $\pmb{u}_B(t)$. We extend the original BCT to allow negative inputs. The resulting piecewise-defined function is defined as: 
\begin{align}
    G(x) = 
\begin{cases}
x^{\lambda_b}  / \lambda_b~ & \text{if } x\geq 0,  \\
-(-x)^{\lambda_b}  / \lambda_b~ & \text{if } x < 0, 
\end{cases}
\end{align}
where $\lambda_b$ is a positive coefficient. In fact, $G(x)$ is the odd extension of $B(x)$, and therefore it inherits the same capability of mapping lower-order quantities to the order of $O(1)$. Moreover, $G(x)$ satisfies the following properties:  
\begin{align}
    G(\prod_{i=1}^n x_i) &= \lambda_b^{n-1}\cdot \prod_{i=1}^{n}G(x_i),~ n\in\mathbb{Z}^{+}  \\
    G( x^{\alpha} ) &= \lambda_b^{\alpha-1}\cdot \left[  G(x) \right]^{\alpha},~ \alpha\in \mathbb{R}^{+}
\end{align}
Applying $G(x)$ to the label $\pmb{u}_B(t)$, the $i$-th entry of $G(\pmb{u}_B(t))$ can be expressed as:
\begin{align}
    G\left( \frac{\mathrm{d} B(Y_{i}) }{ \mathrm{d} t }  \right) 
    &= \lambda_b^{\alpha+2} \cdot  G(C) \cdot G\left(B'(Y_i)\right) 
    \cdot  G\left(k_f\right) 
    \cdot \prod_{k=1}^{N_s} \left[ G(Y_k) \right]^{\alpha_k},
\end{align}
where $\alpha = \sum_{k=1}^{N_s} \alpha_k$. In Eq.~(\ref{eq:bct_ode}), the multiscale terms rate constant $k_f$ and the species concentrations $Y_k$ require further rescaling. After applying $G(x)$, we can tune the parameter $\lambda_b$ such that $G(k_f)$ and $G(Y_k)$ both remain near the order of $\mathcal{O}(1)$. For instance, we assume $Y_1=0.1$, $Y_2=10^{-10}$ (neglecting mass summation conservation for simplicity), and $k_f=10^{9}$. 
With $\lambda_b=0.1$, we obtain $G(k_f)\approx 79$, $G(Y_1)\approx7.9$ and $G(Y_2)=1$.
It is evident that we have successfully transformed quantities originally spanning approximately 20 orders of magnitude to similar orders of magnitude. This transformation effectively mitigates the multiscale effect to a significant extent.

\section{Numerical Results}
\label{sec:result}

In this section, we will conduct comprehensive numerical experiments to systematically evaluate the accuracy, stability, and generalization capabilities of our novel framework. 
We select three representative types of stiff ODE systems of reaction kinetics: 
(1) the well-known 3-species Robertson problem; 
(2) a 21-species methane reaction kinetics model (DRM19\footnote{A. Kazakov, M. Frenklach, Reduced Reaction Sets based on GRI-Mech 1.2, \url{http://www.me.berkeley.edu/drm/}.}); 
and (3) a widely used 13-isotope nuclear reaction network \citep{zhang2025AJ,rivas_impact_2022}.  
For these three types of problems, the DNN-predicted time steps are set to $\Delta t=7\times10^{-7}$ s, $\Delta t=10^{-6}$ s and $\Delta t=10^{-7}$ s, respectively, with the step sizes determined by both stability constraints and practical requirements. 
For all test cases, the BCT hyper-parameter $\lambda_{a}=0.1$ is adopted following Zhang et al. \cite{tz2021}, while the GBCT hyper-parameter is empirically chosen as $\lambda_b=0.5$.

Our benchmark experiments cover a series of increasingly challenging scenarios. To evaluate the basic accuracy and long-term stability of \gbctnet{}, we first examine the pure kinetic evolution of both \bctnet{} and \gbctnet{} on the DRM19 and nuclear reaction models (Sec.~\ref{sec:pure_dynamics}). 
To further assess model performance in cases where reaction kinetics are coupled with complex flow fields, we compare the two networks on the two-dimensional Robertson diffusion problem and a turbulent reaction–diffusion system (Sec.~\ref{sec:chemical_flow}). 
In addition, we investigate their capabilities for prediction in complicated astrophysical reactive flows characterized by extreme stiffness (Sec.~\ref{sec:nuclear_flow}). This comprehensive selection ensures that our benchmarks encompass both canonical stiff systems and complex, practically relevant cases involving coupled transport phenomena. The code is available at the GitHub repository\footnote{\url{ https://github.com/Seauagain/GBCT}} for reproduction.

\subsection{Reaction Dynamics Simulations}
\label{sec:pure_dynamics}

\subsubsection{  Case 1: methane/air reaction dynamics}

\begin{figure}[htbp]
    \centering
    \includegraphics[width=1\linewidth]{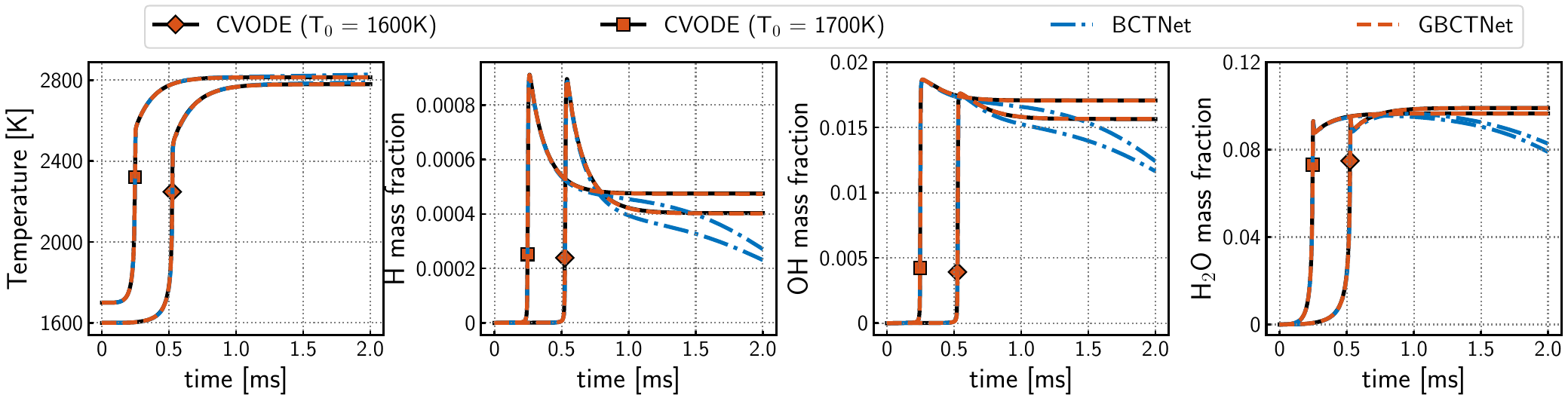}
    \caption{
     Temporal evolution of temperature and representative radicals \ce{H}, \ce{OH},\ce{H2O} of two test cases with initial temperature $T_0 = 1600$ K and  $T_0 = 1700$ K, respectively. The initial pressure is $p_0 = 1$ atm and equivalence ratio is $\phi_0 = 1.0$. The black, red, and blue lines represent the simulation results from the CVODE, \gbctnet{} and \bctnet{}, respectively.
        }
    \label{fig:drm19_0d}
\end{figure}

We begin with validating the prediction capability of the model for chemical kinetics. 
Specifically, we investigate the auto-ignition process of a methane/air mixture within a constant-pressure reactor, which is described by a set of 21-dimensional stiff ordinary differential equations (ODEs). Two test scenarios are considered, both at an initial pressure of $p_0 = 1$ atm and an equivalence ratio of $\phi_0 = 1.0$, but differing in initial temperatures of $T_0 = 1600$ K and $T_0 = 1700$ K, respectively. The reaction dynamics are simulated over 2000 discrete time steps. During this process, the deep neural network (DNN) models iteratively predict the mass fractions of chemical species at each subsequent time step, while simultaneously updating the temperature based on energy conservation principles. As illustrated in Fig. \ref{fig:drm19_0d}, both DNN models successfully capture the temporal evolution of temperature. Nevertheless, the \bctnet{} model exhibits significant inaccuracies in predicting the mass fractions of radicals such as \ce{H}, \ce{OH}, and \ce{H2O}, primarily due to the accumulation of errors over time.

\begin{figure}[htbp]
    \centering
    \includegraphics[width=1\linewidth]{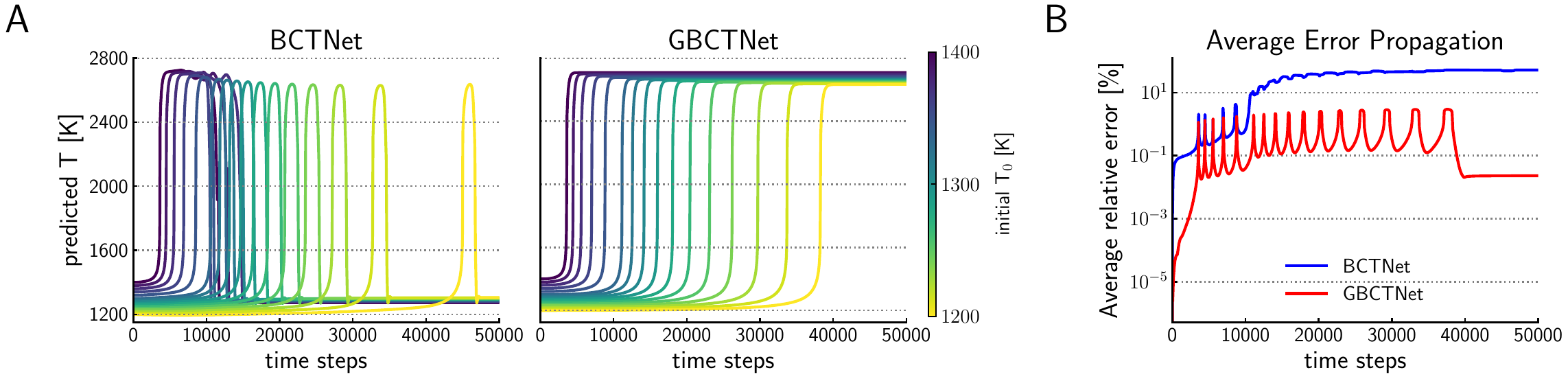}
    \caption{ A) Comparison of the long-term temperature evolution between \bctnet{} (left) and \gbctnet{} (right) across varying initial temperatures $T_0 \in [1200, 1400]$ K. (B) The average relative error propagation in temperature. The initial pressure is set to $p_0 = 1$ atm, with an equivalence ratio of $\phi_0 = 1.0$. The simulation is run for 50,000 steps, corresponding to 50 ms since the time step  $\Delta t = 10^{-6}$ s used for DNN prediction. }
    \label{fig:drm19_longterm}
\end{figure}

To further assess the model's accuracy and stability over long-term evolution, we uniformly sample 16 initial points from the temperature range $T_0 \in [1200, 1400]$ K. 
The corresponding temporal trajectories required a significantly prolonged integration period to reach steady state, necessitating extended simulation durations to fully capture the system's evolutionary characteristics. 
Fig. \ref{fig:drm19_longterm} (A) compares the long-term temperature evolution between \gbctnet{} and \bctnet{}, where curves of identical colors share the same initial temperature, with fixed initial pressure $p_0 = 1$ atm and equivalence ratio $\phi_0 = 1.0$. 
The simulation spans 50,000 time steps (equivalent to 50 ms) to reach steady state. 
The results demonstrate that \gbctnet{} exhibits remarkable long-term stability, consistently maintaining accurate temperature predictions even after 50,000 iterative steps.
Quantitatively, we define the relative error of equilibrium temperature as $\vert T_{\text{eq}} - \tilde{T}_{\text{eq}} \vert / T_{\text{eq}} $, where $T_{\text{eq}}$ denotes the ground truth obtained via CVODE and  $\tilde{T}_{\text{eq}}$ represents the model's prediction. 
As shown in Fig. \ref{fig:drm19_longterm} (B), 
\gbctnet{} achieves a mean relative error for equilibrium temperature of about 0.1\%, compared to approximately 50\% for \bctnet{}, highlighting \gbctnet{}'s superior performance.

\subsubsection{ Case 2: astrophysical nuclear reaction}

\begin{figure}[htbp]
    \centering
    \includegraphics[width=1\linewidth]{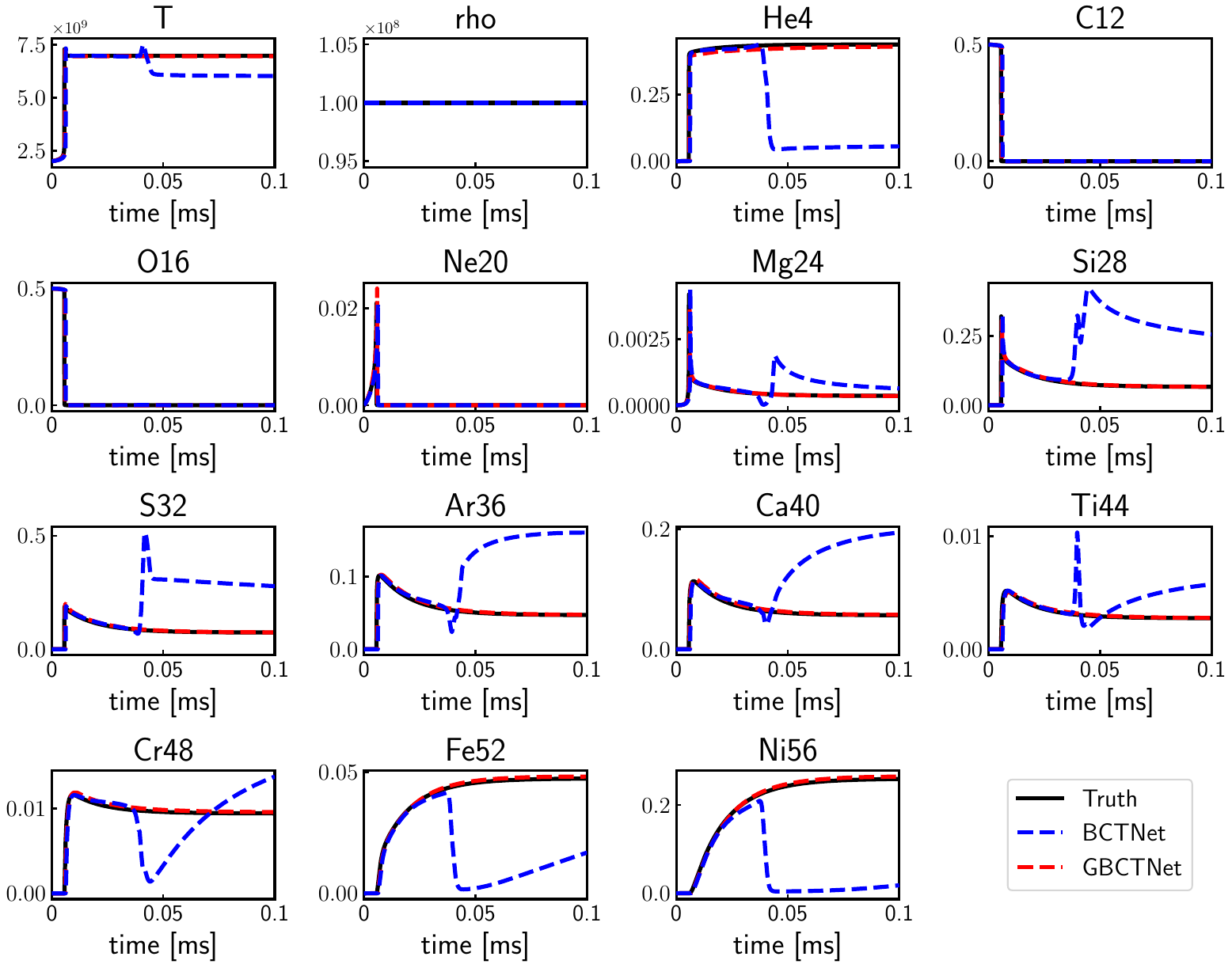}
    \caption{The evolution of $^{12}$C-$^{16}$O nuclear reaction dynamics. The initial temperature is $T_0 = 2\times10^9$ K and density $\rho_0 = 1\times 10^8$ g/cm$^3$ with simulation time of 0.1 ms. The black solid line denotes simulations obtained by CVODE. The red and blue dashed line denote the predicted results of \gbctnet{} and \bctnet{}, respectively. 
    }
    \label{fig:0d_nuclear}
\end{figure}

Another representative example involves astrophysical nuclear reactions. We select a 13-isotope reaction kinetic model that is commonly employed in combustion simulations of helium shells and carbon-oxygen cores. These reactions typically take place under extreme conditions, with temperatures on the order of $\sim \mathcal{O}(10^{9})$ K and densities around $\sim \mathcal{O}(10^{8})$ g/cm$^3$, which require extremely small temporal integration steps due to the rapid reaction rates.

Our test case initializes with thermodynamic parameters $T_0 = 2 \times 10^9$ K and $\rho_0 = 1 \times 10^8$ g/cm$^3$, and the nuclide mass fraction consists of 50\% $^{12}$C and 50\% $^{16}$O. Numerical integration is carried out using the CVODE solver over a duration of 0.1 ms. The DNN model simulations employ a time step of $\Delta t = 10^{-7}$ s. As illustrated in Fig. \ref{fig:0d_nuclear}, \gbctnet{} accurately captures the temporal evolution of temperature and isotopic mass fractions. In contrast, \bctnet{} fails to correctly predict the trajectories within the reaction equilibrium regions.

\subsection{Chemical Diffusion-Reaction Dynamics Simulations}
\label{sec:chemical_flow}

In diffusion-reaction systems, the reaction dynamics play a pivotal role in influencing heat release and pressure fluctuations, thereby altering the flow characteristics. 
These changes in flow dynamics, in turn, have a profound impact on temperature and the spatiotemporal distribution of reactants. Such intricate interactions between flow and reaction necessitate a precise resolution of reaction dynamics.
On the other hand, decoupled reaction kinetics solutions on separate grids can be independently computed, enabling large-scale, grid-level parallel computations. As a result, diffusion-reaction systems, particularly under complex flow configurations, provide an ideal framework for validating both the accuracy and acceleration potential of surrogate models for reaction kinetics.
In this subsection, we integrate the DNN model into the solver for diffusion-reaction systems, replacing traditional ODE numerical integrators.

\subsubsection{Case 3: two-dimensional Robertson diffusion problem}
\label{sec:case3_2d_robert}

The Robertson problem (abbreviated as ROBER) \cite{robert1966} is a typical stiff dynamical system. We consider the Robertson problem coupled with diffusion, whose governing equations are:
\begin{align}
& \frac{\partial y_1}{\partial t} = D_1\nabla^2 y_1 -k_1y_1 +k_3y_2y_3,\label{eq:y1} \\
& \frac{\partial y_2}{\partial t} = D_2\nabla^2 y_2 + k_1y_1 - k_2y_2^2 - k_3 y_2 y_3, \label{eq:y2} \\
& \frac{\partial y_3}{\partial t} =D_3\nabla^2 y_3 + k_2y_2^2. \label{eq:y3}
\end{align}
The ROBER ODEs describe the chemical reaction process involving 3 species, where the reaction rates exhibit vastly different orders of magnitude. In these equations, $y_i$ represents the concentrations of the 3 species, while the reaction rate constants are $k_1 = 4\times 10^{-2}$, $k_2=3\times10^{7}$, and $k_3 = 10^{4}$, respectively. The diffusion coefficients are set as $D_1=100$, $D_2= 0.5$, and $D_3 = 2$. The computational domain is set as $[-1,1]\times[-1,1]$. More details about the configuration are provided in the Append. \ref{sec:robert_appd}.

\begin{figure}[htbp]
	\centering
	\includegraphics[width=1.0\linewidth]{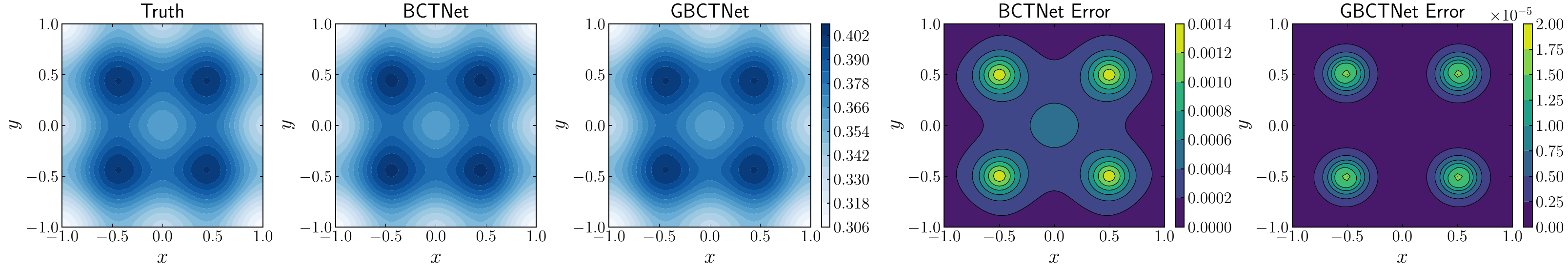}
	\caption{The species concentration $y_1(x,y)$ and the absolute error of BCTNet and GBCTNet for the Robertson diffusion problem after 1000 simulation steps with a time step of $\Delta t=7\times10^{-7}$ s.}
	\label{fig:robert_2d}
\end{figure}
For the numerical solution, the Strang splitting method \cite{Strang1968} is employed to decouple the diffusion and reaction operators. BCTNet and GBCTNet with a time step $\Delta t = 7\times 10^{-7}$ s are trained to replace the integral of stiff reaction dynamics, respectively. The simulation is executed for 1000 time steps, using results obtained with the Runge-Kutta method as the reference solution. As shown in the left panel of Fig.~\ref{fig:robert_2d}, the component concentration $y_1(x,y)$ predicted by BCTNet and GBCTNet closely matches the true values. The right panel of Fig.~\ref{fig:robert_2d} further illustrates the absolute error distribution of the two networks. Compared to BCTNet, GBCTNet achieves an absolute error on the order of approximately $\mathcal{O}(10^{-5})$, which is two orders of magnitude lower than the former, demonstrating its high accuracy and stability in continuous predictions.

\subsubsection{Case 4: two-dimensional turbulent ignition}

After demonstrating \gbctnet{}'s  capability in predicting qualities of interest in ROBER diffusion problem, we further investigate its performance in a more challenging scenario: two-dimensional turbulent chemically reaction-diffusion (i.e. combustion). The details of simulation code and governing equations can be found in Append.~\ref{sec:ebi_equation_appd}.

In this case, the computational domain is set as a rectangular region of 
1.5 cm $\times$ 1.5 cm. The grids for discretization is N=512 $\times$ 512. Initially, the computational area is filled with a mixture of \ce{CH4}, \ce{O2}, and \ce{N2} in a volume ratio of $1:1:3.76$. The initial temperature $T_0$ is 300 K, and the initial pressure $p_0$ is 1 atm. To generate an initial isotropic turbulent field, we employ the Passot-Pouquet isotropic kinetic energy spectrum, which can be expressed as
\begin{align}
&E(k)=16 \sqrt{\frac{2}{\pi}} \frac{u_{rms}^2}{k_e}\left(\frac{k}{k_e}\right)^4 \exp \left(-2\left(\frac{k}{k_e}\right)^2\right) , 
\end{align}
where $E(k)$ is the energy spectrum, $k$ is the frequency, the average velocity $u_{rms}=3$ m/s, and $k_e=418.67$. Given the energy spectrum, we first generate a velocity field in spectral space which satisfies the continuity equation, and then transform it into physical space via the inverse Fourier transform \cite{saad2017}. The ignition source, positioned at $(0.75\,\mathrm{cm}, 0.75\,\mathrm{cm})$, features a radius of $0.02$ cm and energy density of $1\times10^{11} \mathrm{W/m^3}$, activated for $0.2$ ms duration.

\begin{figure}[htbp]
    \centering
    \includegraphics[width=1.0\linewidth]{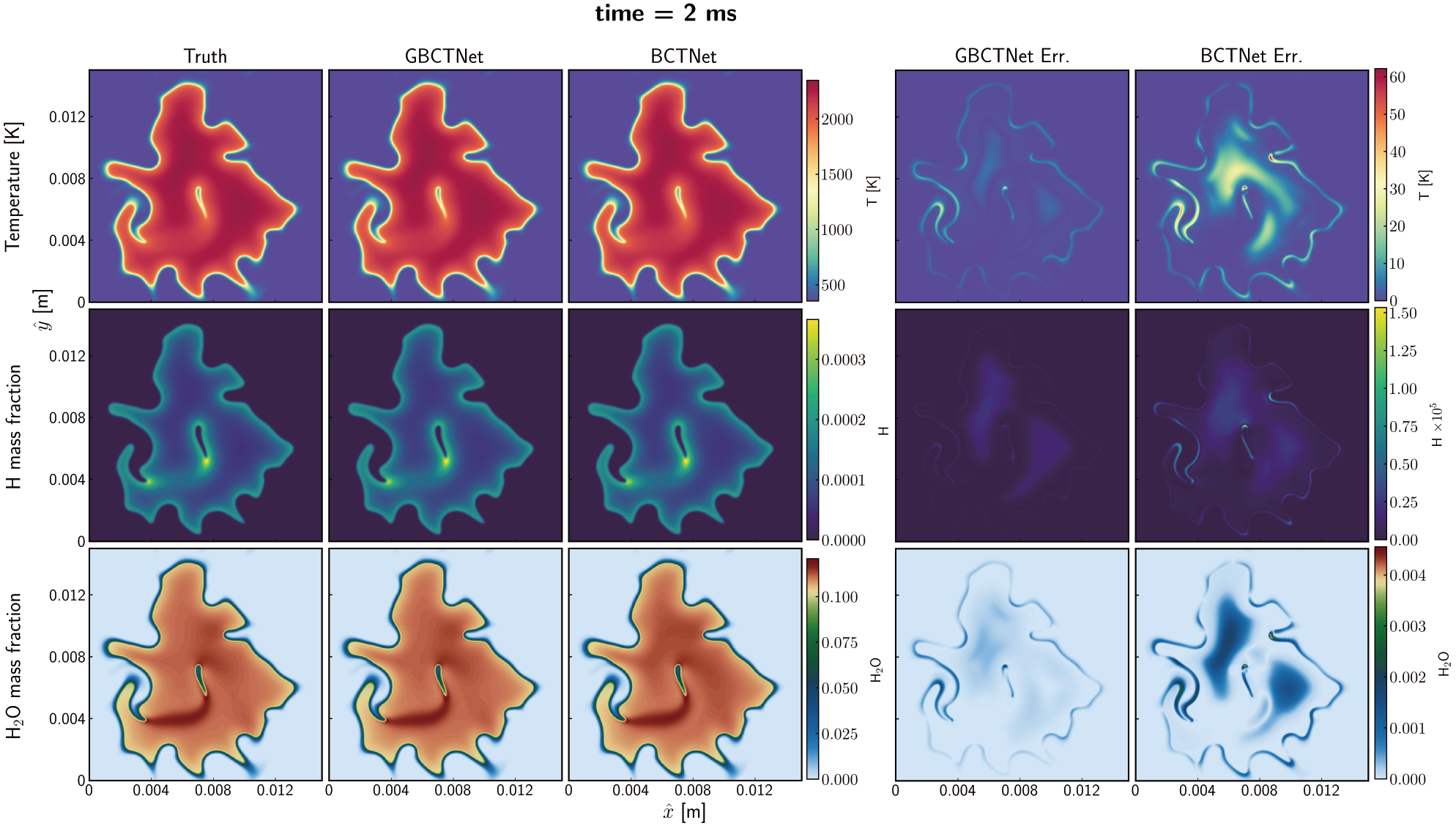}
    \caption{Two-dimensional turbulent case: comparison of temperature and radical mass fraction snapshots (left) and absolute error comparison (right) at t = 2 ms.}
    \label{fig:drm19_2d}
\end{figure}

The DNN models are integrated into grid cells where temperatures exceed 800 K to predict species variations caused by chemical reactions. 
Fig. \ref{fig:drm19_2d} presents snapshots of temperature and radical mass fractions predicted by CVODE and the two DNN models at time = 2 ms. Both \bctnet{} and \gbctnet{} demonstrate reasonable accuracy in capturing the contour profiles. A closer examination of the absolute error distributions (Fig. \ref{fig:drm19_2d}, right panel) reveals that \gbctnet{} achieves significantly smaller absolute errors compared to \bctnet{}. Specifically, the maximum absolute error in temperature predicted by \gbctnet{} is 30 K, which is approximately half of the maximum error observed in \bctnet{} predictions. The primary errors in \bctnet{} are concentrated in high-temperature regions near the center of the computational domain (temperatures > 2000 K), particularly for temperature and the mass fraction of the product \ce{H2O}. This suggests that \bctnet{} struggles with accuracy in burned regions, whereas \gbctnet{} effectively reduces errors in these areas. These findings align with our earlier conclusions.

\begin{figure}[htbp]
    \centering
    \includegraphics[width=1.0\linewidth]{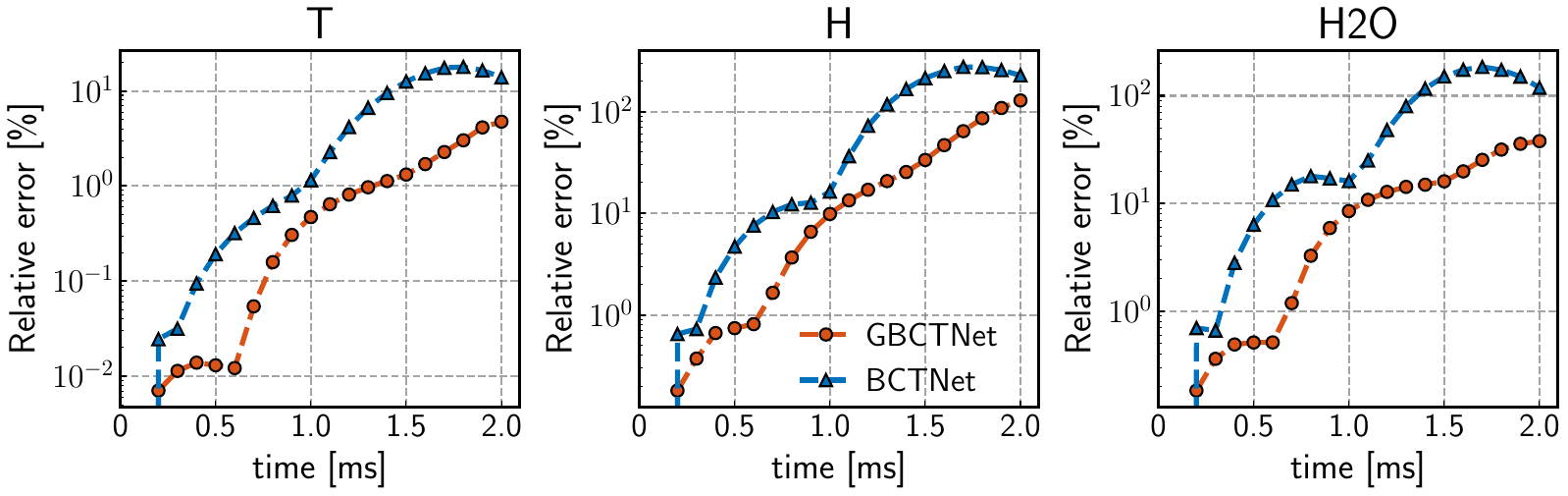}
    \caption{Two-dimensional turbulent case: the domain-average relative error propagation over time of the \gbctnet{} (red circles) and \bctnet{} (blue triangles) prediction.}
    \label{fig:drm19_2d_error}
\end{figure}

To quantitatively compare the temporal error propagation of the two DNN models from 0 ms to 2 ms, we systematically calculated the domain-averaged relative errors between neural network predictions and CVODE reference solutions at 0.1 ms intervals (Fig. \ref{fig:drm19_2d_error}). Both models demonstrate comparable accuracy, with temperature predictions generally maintaining relative errors below 20\%. However, \gbctnet{} exhibits consistently lower errors than \bctnet{} throughout the simulation period. At t = 2 ms, \gbctnet{} achieves a temperature relative error of 4.8\%, outperforming \bctnet{}'s 17.8\%.
For the chemically active and error-sensitive radical species \ce{H}, \gbctnet{} attains a maximum relative error of 100.2\%, substantially lower than \bctnet{}'s 274\%. Similarly, \gbctnet{} reduces the maximum relative error in \ce{H2O} mass fraction predictions from 185\% (\bctnet{}) to 37\%. These results quantitatively confirm \gbctnet{}'s enhanced robustness in handling both stable products and reactive intermediates across the turbulence-reaction process.

\subsection{Nuclear Diffusion-Reaction Dynamics Simulations}
\label{sec:nuclear_flow}

To display the versatility of our method, we extend the application of our model to simulate astrophysical reacting flows. 
While similar to chemical combustion, nuclear reactions occur under far more stringent conditions, characterized by immense energy release and significantly shorter characteristic timescales compared to chemical reactions.
Therefore, accurately capturing reacting fronts necessitates precise modeling of nuclear reactions, especially when coupled with hydrodynamic processes.
In this subsection, we present and analyze two hydrodynamic simulation cases: one-dimensional unsteady flame propagation and two-dimensional flame simulation.

The governing equations for reacting flow simulations can be expressed as:
    \begin{align}
    &\frac{\partial (\rho Y_{k})}{\partial t} + \nabla \cdot (\rho Y_{k} \bm{v}) = \rho 
    \dot{\omega}_{{k}},  \label{eq:mass} \\
   &\frac{\partial(\rho\bm{v})}{\partial t} + \nabla \cdot(\rho\bm{v}\bm{v}+p\bm{I}) = 0 ,  \label{eq:momentum} \\
   &\frac{\partial(\rho\epsilon)}{\partial t} + \nabla \cdot[\bm{v}(\rho\epsilon+p)] =  \rho \dot{S}. \label{eq:energy}
    \end{align}
Here $\rho$, $p$, $\bm{v}$ are the fluid mass density, gas pressure and velocity, respectively, $Y_{k}$ is the mass fraction of species $k$ with associated production rate $\dot{\omega}_{{k}}$ and $\bm{I}$ is the identity tensor. In Eq.(\ref{eq:energy}), $\epsilon$ denotes the total specific energy and $\dot{S}$ is the nuclear energy generation rate per unit mass.

Typically, nuclear reactions and advection are often separated using operator splitting techniques, such as Strang splitting \cite{Strang1968}, allowing them to be computed independently. 
The governing ordinary differential equations (ODEs) for nuclear reactions are expressed as:
\begin{align}
& \frac{\mathrm{d} \bm{Y} }{\mathrm{d} t} = \bm{\dot{\omega}} (T, \rho, \bm{Y}), \\
& \frac{\mathrm{d} (\rho\epsilon)}{\mathrm{d} t}  =  \rho e_{nu},
\end{align}
where $\bm{Y}$ represents the vector of species mass fractions, and the production rates $\bm{\dot{\omega}}$ are typically dependent on temperature and density. 
The temperature $T$ is updated via the equation of state (EOS). For a gamma-law EOS:
\begin{align}
p = \rho \epsilon_{int} (\gamma-1).
\end{align}
Here, the ideal gas law is applied with an adiabatic coefficient of $\gamma=5/3$. The specific internal energy $\epsilon_{int}$ is determined by subtracting the specific kinetic energy from the total specific energy $\epsilon$.
In conventional computational workflows, the time step for advancing nuclear reaction networks is considerably smaller compared to that for hydrodynamic processes (i.e., $\Delta t_{\text{react}} \ll \Delta t_{\text{flow}}$), resulting in significant computational costs. This study develops surrogate models to overcome stiffness constraints, thereby efficiently replacing the time-intensive direct integration (DI).

\subsubsection{Case 5: one-dimensional unsteady nuclear flame front propagation}

\begin{figure}[htbp]
    \centering
    \includegraphics[width=1.0\linewidth]{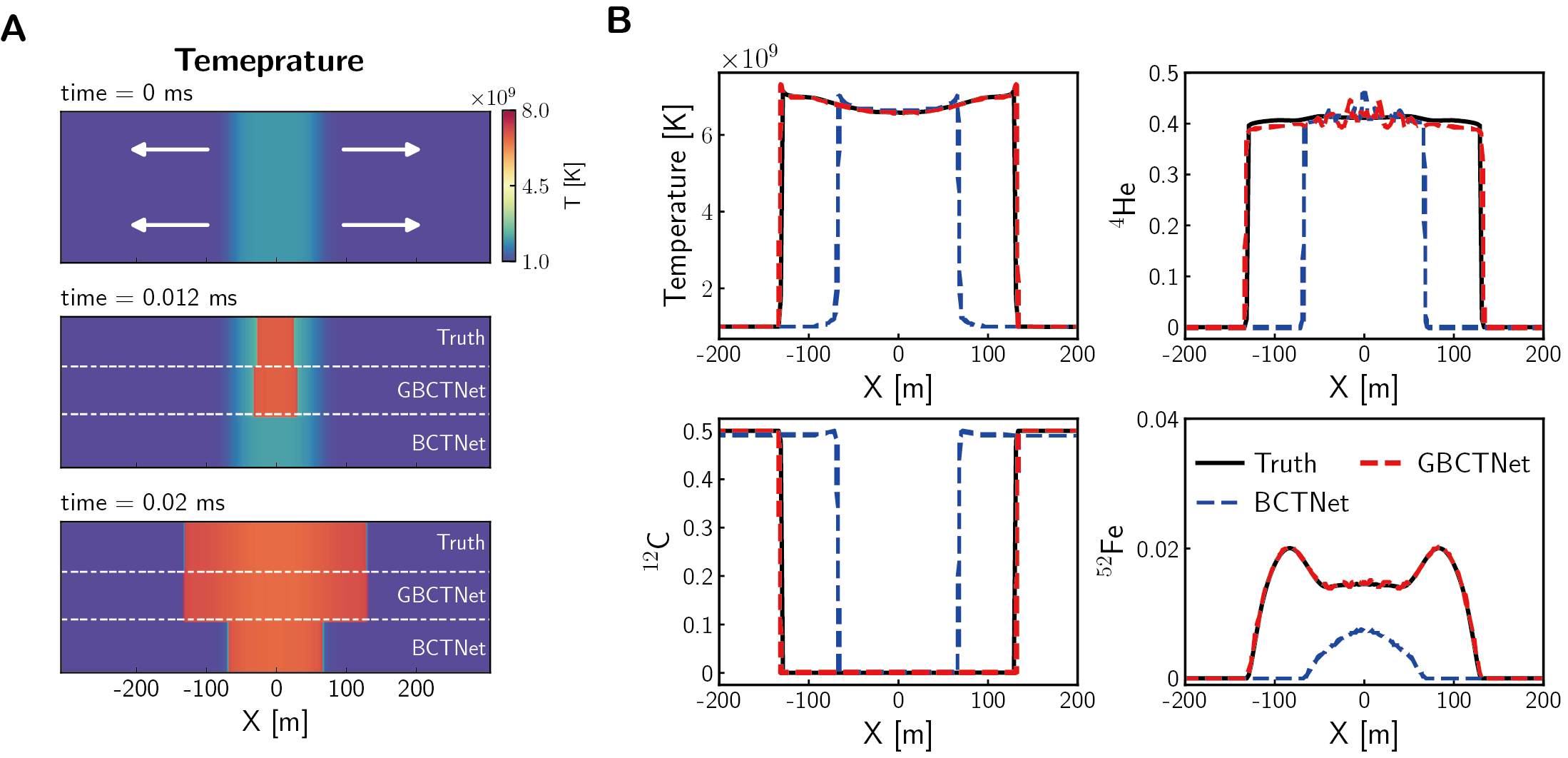}
    \caption{One-dimensional unsteady nuclear flame simulations. (A) Temperature snapshots obtained by \bctnet{}, \gbctnet{} and ground truth at $t=0$, 0.012, 0.02 ms, respectively. (B) Temperature and isotope concentration distributions predicted by \gbctnet{} and \bctnet{} models at $t=2$ ms.
 }
    \label{fig:nuclear_1d}
\end{figure}

We first perform a one-dimensional unsteady flow simulation to describe the internal deflagration process in white dwarfs. Unlike the scenario discussed in Sec. \ref{sec:case3_2d_robert}, this benchmark necessitates modeling the propagation of nuclear reaction fronts and accurately predicting nuclear burning regions. 
A one-dimensional computational domain of 614.4 m is established, with outflow boundary conditions applied at both ends.
The Cartesian coordinate system is centered at the midpoint, and the domain is discretized into 256 grid cells ($\Delta x = 15.37$ m). 
Initially, the domain consists of an isotopic composition of approximately 50\% $^{12}C$ and 50\% $^{16}O$, with a uniform density of $\rho = 1\times 10^7$ g/cm$^3$. 
The upper row of Fig. \ref{fig:nuclear_1d}(A) illustrates the initial configuration. 
The temperature field is defined as $T(x) = \alpha(x) \cdot T_0$, where $T_0 = 2\times10^8$ K and the scaling factor $\alpha(x) = 1.5 + 0.5\tanh(4 - |x|/\Delta x)$. 
This setup results in the maximum temperature at the origin, with a gradual decay toward both boundaries. Ignition begins near the coordinate origin, forming a nuclear flame front that propagates in both directions.

The middle and bottom rows in Fig.~\ref{fig:nuclear_1d}(A) depict temperature distributions at $t = 0.012$ ms and $t = 0.02$ ms, respectively, obtained from \gbctnet{}, \bctnet{}, and the numerical solver. At $t = 0.012$ ms, the flame front (red region) predicted by \gbctnet{} shows initial propagation with excellent agreement to the reference solution, while \bctnet{} fails to initiate combustion. At $t = 0.02$ ms, \gbctnet{} continues to align closely with the reference solution, whereas \bctnet{} demonstrates significantly slower flame front propagation.

Fig.~\ref{fig:nuclear_1d}(B) presents quantitative comparisons of temperature and three representative isotopic concentrations at $t = 0.02$ ms. \bctnet{} predictions show significant deviations from the ground truth, with a mean relative temperature error of 17.8\%. In contrast, \gbctnet{} exhibits superior predictive accuracy in terms of flame front position, peak reaction temperature, and isotopic concentrations. The relative temperature error for \gbctnet{} is measured at 2.6\%, underscoring its excellence in predicting nuclear dynamics.

\subsubsection{Case 6: Two-dimensional wedge-shaped nuclear flame front propagation}

\begin{figure}[htbp]
    \centering
    \includegraphics[width=1.0\linewidth]{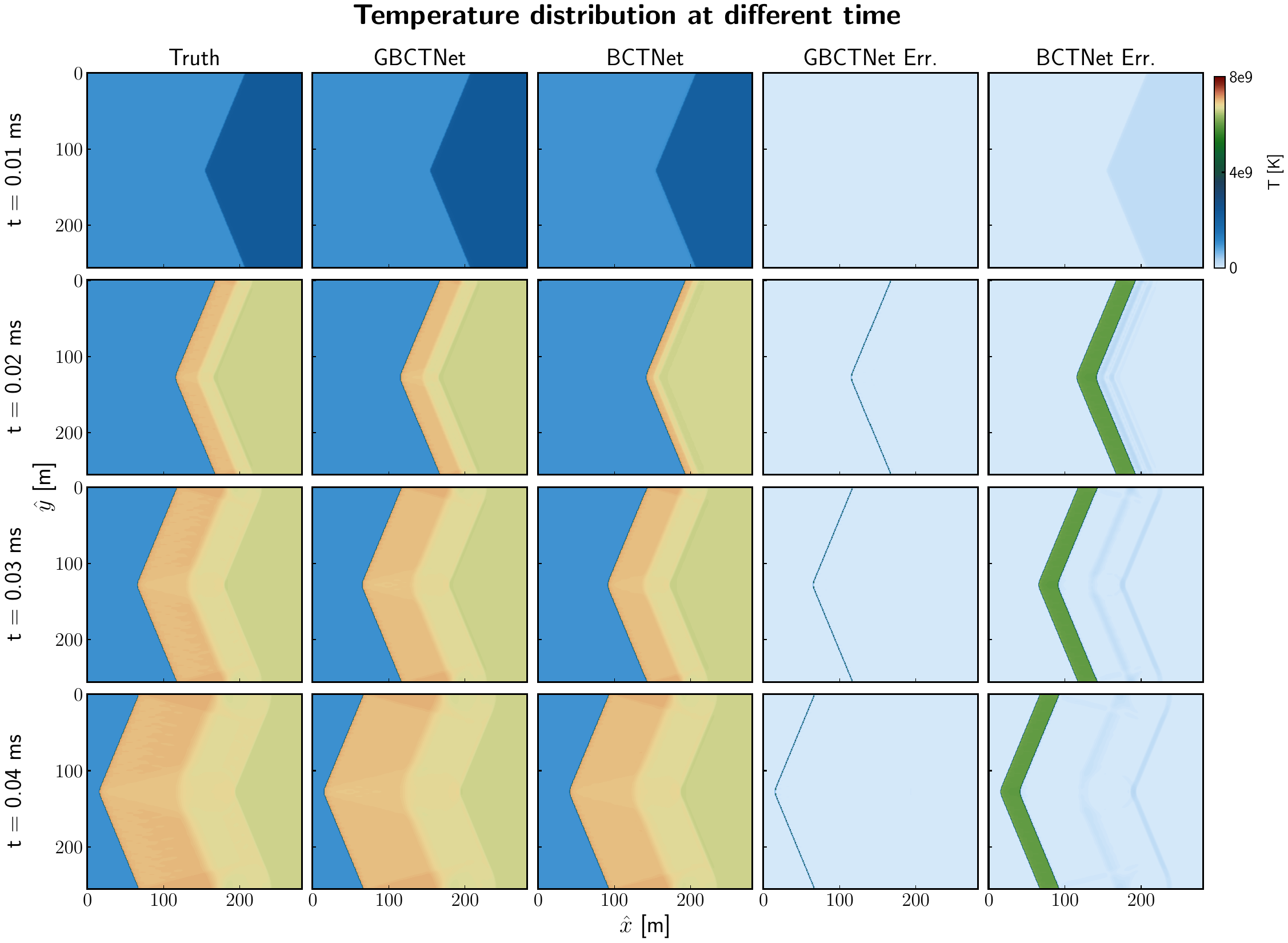}
    \caption{Two-dimensional nuclear flame. Temperature snapshots predicted by two DNN models (\bctnet{} and \gbctnet{}) at $t = $ 0.01, 0.02, 0.03, and 0.04 ms, along with the absolute error compared with the ground truth.}
    \label{fig:nu2d_contour}
\end{figure}

\begin{figure}[htbp]
    \centering
    \includegraphics[width=1.0\linewidth]{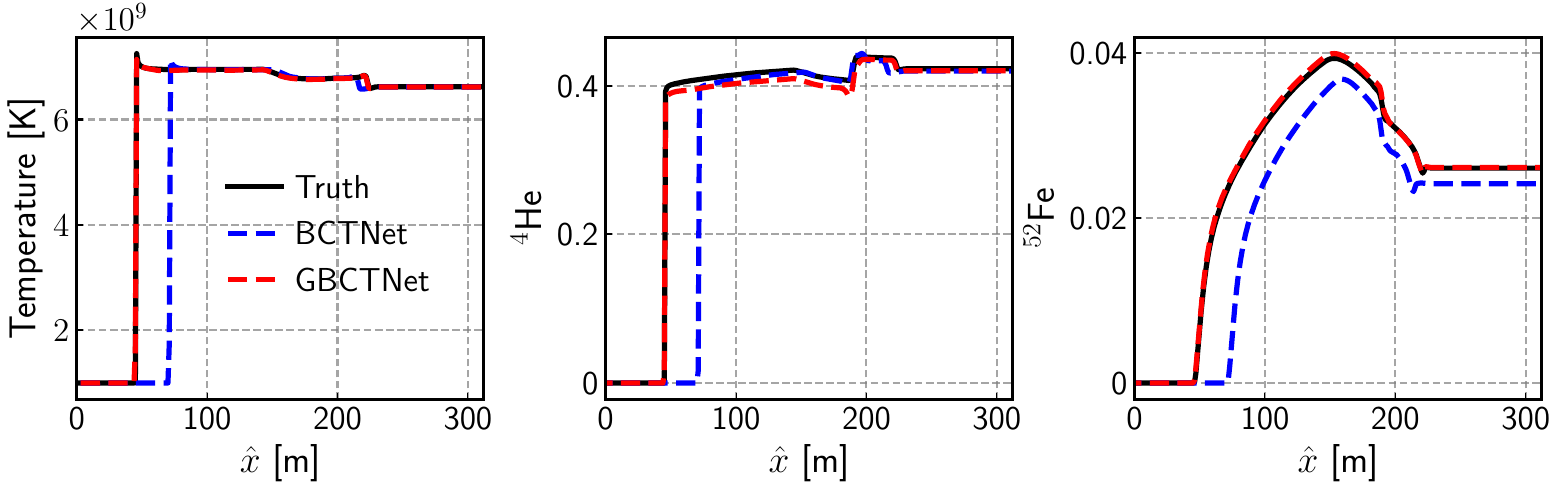}
    \caption{ Temperature and isotope concentration profiles along the horizontal mid-plane ($y = H/2$) at $t = 0.04$ ms.}
    \label{fig:nu2d_line}
\end{figure}

We delve into a more intricate benchmark test involving two-dimensional nuclear flame propagation.
The computational domain is a two-dimensional rectangle measuring 921.6 m $\times$ 614.4 m, discretized into a grid of 384 $\times$ 256 cells, with a resolution of $\Delta x = \Delta y = 2.4$ m. Boundary conditions are defined as outflow/reflective along the x-axis and periodic along the y-axis. Initially, the domain is filled with a homogeneous mixture comprising 50\% $^{12}$C and 50\% $^{16}$O. A V-shaped boundary divides the domain into left and right sub-regions, with the upper and lower edges forming angles of $67.5^\circ$ relative to the x-axis. The right subregion is initialized with $T_1 = 2\times10^9$ K and $\rho_1 = 5\times 10^7$ g/cm$^3$, while the left subregion starts with $T_2 = 1\times 10^8$ K and $\rho_2 = 1\times10^8$ g/cm$^3$. Nuclear combustion begins in the preheated right subregion, generating a wedge-shaped flame front that advances into the unburned left region.

Fig.~\ref{fig:nu2d_contour} illustrates temperature snapshots predicted by two DNN models (\bctnet{} and \gbctnet{}) at $t$ = 0.01, 0.02, 0.03, and 0.04 ms, alongside absolute error maps relative to the ground truth, respectively.
Both models exhibit stable predictions of the flame front during propagation. However, \bctnet{} produces notable absolute errors, characterized by spatially extensive error regions near the flame front. In contrast, \gbctnet{} demonstrates superior accuracy, with errors confined to a narrow band closely aligned with the flame front.

Fig.~\ref{fig:nu2d_line} provides a quantitative comparison of temperature and isotope concentration profiles along the horizontal mid-plane ($y = H/2$) at $t = 0.04$ ms. The quantitative analysis aligns with the snapshot observations: \bctnet{} predictions (blue dashed lines) show significant deviations at the flame front, whereas \gbctnet{} predictions (red dashed lines) accurately replicate the physical parameter distributions.

\begin{figure}[htbp]
    \centering
    \includegraphics[width=1.0\linewidth]{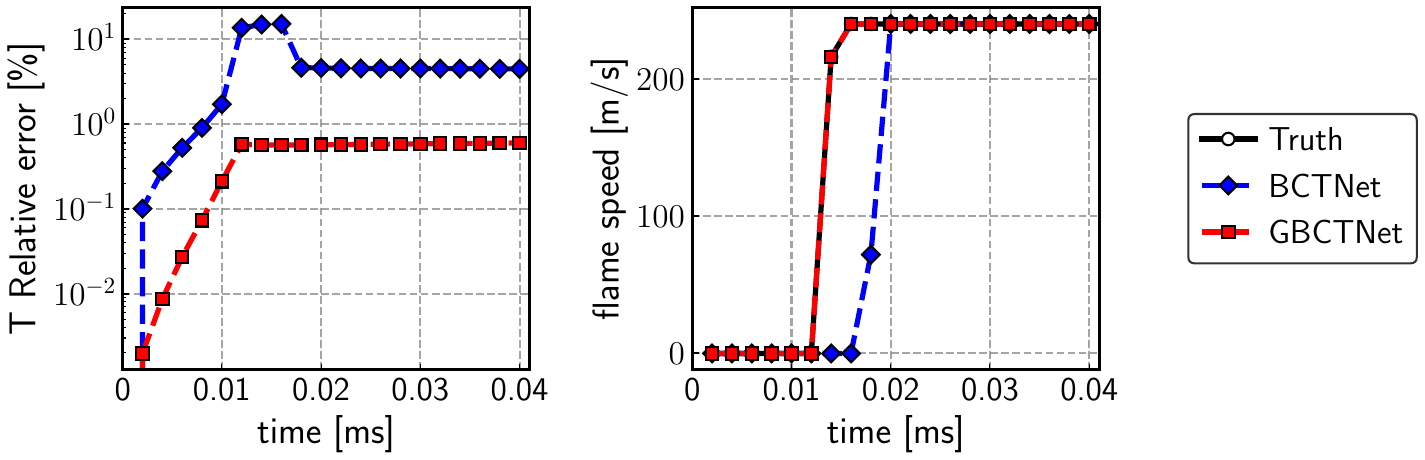}
    \caption{Two-dimensional nuclear reaction flame. The domain-average relative errors and flame velocity obtained by \bctnet{} and \gbctnet{} at different time.  }
    \label{fig:nu2d_error}
\end{figure}

Fig.~\ref{fig:nu2d_error} compares domain-average relative errors and flame propagation velocities. \gbctnet{} achieves temperature prediction errors that are an order of magnitude lower than those of \bctnet{}. The ground truth flame front initiates propagation at 0.01 ms, rapidly accelerating to a velocity of 200 m/s. Although \bctnet{} eventually matches the final propagation speed, it exhibits delayed ignition timing, which accounts for its elevated errors. In contrast, \gbctnet{} demonstrates excellent agreement with both the temporal evolution and spatial progression of the flame front throughout the simulation.

\section{Analysis}
\label{sec:analysis}

\subsection{Impact of GBCT on training efficiency}
\label{sec:training}

This section investigates the impact of GBCT on the training efficiency of neural networks and explores the underlying mechanism through the lens of frequency.
For clarity, we employ the chemical kinetics described in Sec.~\ref{sec:pure_dynamics} as a representative test example.
To minimize random fluctuations, both \bctnet{} and \gbctnet{} are trained for five trials under identical settings with different random seeds.

\begin{figure}[htb]
    \centering
    \includegraphics[width=1\linewidth]{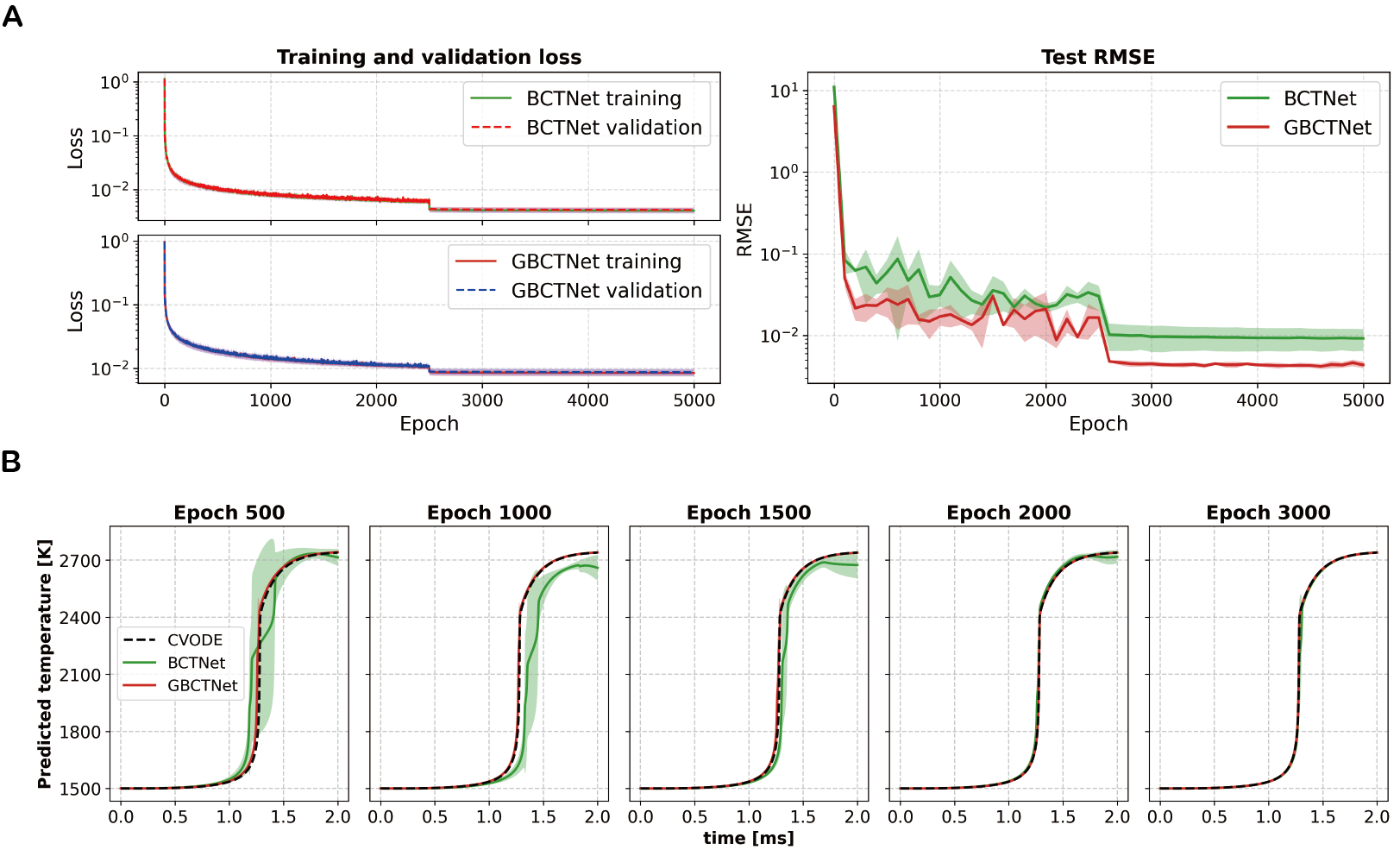}
    \caption{ The influence of GBCT on neural network training. (A) Left panel: training and validation losses for both networks; Right panel: RMSE on the unseen test set. (B)  Continuous temporal prediction of the neural network models at different training epochs. Testing is conducted under initial conditions of $T_0$ = 1500 K, $p_0$ = 1 atm, equivalence ratio $\phi_0$ = 1.0, with 2000 testing steps. Each setup of training is conducted for 5 trials with different random seeds.
    }
    \label{fig:train_analysis}
\end{figure}

Fig.~\ref{fig:train_analysis}(A) presents the training loss, validation loss, and RMSE on an unseen test set across training epochs for both networks. The test set consists of samples drawn from thousands of trajectories of temporal evolution.
A two-phase training strategy is adopted: at epoch 2500, the learning rate is reduced while the batch size is increased to accelerate convergence.
As illustrated in the left panel of Fig.~\ref{fig:train_analysis}(A), both training and validation losses decrease significantly following this adjustment and remain stable throughout the subsequent 2500 epochs.
The final training losses for \bctnet{} and \gbctnet{} reach the same order of magnitude $\mathcal{O}(10^{-2})$. 
The right panel of Fig.~\ref{fig:train_analysis}(A) shows that \gbctnet{} consistently achieves lower single-step prediction RMSE on the test set compared to \bctnet{} throughout the training process.
Furthermore, Fig.~\ref{fig:train_analysis}(B) compares the multi-step iterative prediction performance of both models at several training checkpoints (epochs 500, 1000, 1500, 2000, and 3000).
Even during the early training stages (e.g., around epoch 500), its predictions closely match the ground truth obtained from the CVODE solver, whereas the baseline model requires approximately 3000 epochs to reach a comparable level of predictive performance.
These findings confirm that GBCT enhances both training efficiency and predictive performance, including accuracy and robustness.

\subsection{Frequency analysis}
\label{sec:freq}

The remarkable effectiveness of GBCT in both training and prediction inspires further investigation into its underlying mechanisms.
Due to the frequency principle \cite{Xu2020CiCP, xu2025CAMC, xu2019NIP} or spectral bias \cite{rahaman19a}, neural networks tend to prioritize fitting low-frequency components of the target function, while convergence slows when high-frequency components are predominant.
 Thus, we analyze the frequency shift in the training data to elucidate how GBCT influences the training process.


\subsubsection{Two frequency analysis indicators for high-dimensional functions }

Frequency analysis of high-dimensional functions poses a significant challenge. To tackle this, we begin by introducing two effective metrics. Without loss of generality, the Fourier transform of the function $\boldsymbol{y}(\boldsymbol{x})\in L^2(\mathbb{R}^d)$ is defined as:
\begin{align}
    \hat{ \boldsymbol{y} }( \boldsymbol{k} ) = \int_{\mathbb{R}^d} \boldsymbol{y}( \boldsymbol{x} ) \mathrm{e} ^{-i\cdot 2\pi \boldsymbol{k} \boldsymbol{x}} \mathrm{d} \boldsymbol{x}.
\end{align}

We aim to extract the low-frequency components of $\boldsymbol{y}(\boldsymbol{x})$, but the Fourier transform of high-dimensional functions suffers from the curse of dimensionality. Therefore, similar to previous work \cite{Xu2020CiCP, xu2021AAAI}, we use a filtering method to address this issue. 
The basic idea is to separate the function into low-frequency part with frequency $|\boldsymbol{k}| \leq k_0$ and high-frequency part with frequency $|\boldsymbol{k}| > k_0$, using the designated threshold $ k_0\in\mathbb{R} $. 
We then define the Low Frequency Ratio (LFR) as:
\begin{align}
    \text{LFR} (k_0) = \frac{  \displaystyle \int_{\mathbb{R}^d} 
 \left\vert \hat{ \boldsymbol{y} }( \boldsymbol{k} ) \cdot \mathbbm{1}_{|\boldsymbol{k}|\leq k_0} \right\vert ^2\mathrm{d} \boldsymbol{k} } { \displaystyle \int_{\mathbb{R}^d} 
 \left\vert \hat{ \boldsymbol{y} }( \boldsymbol{k} )  \right\vert ^2\mathrm{d} \boldsymbol{k} } ,
\end{align}
where $\mathbbm{1}_{|\boldsymbol{k}|\leq k_0}$ denotes the indicator function, i.e., 
\begin{align}
    \mathbbm{1}_{|\boldsymbol{k}|\leq k_0} = 
    \begin{cases} 
        1, & |\boldsymbol{k}| \leq k_0 \\
        0, & |\boldsymbol{k}| > k_0
    \end{cases}
\end{align}
The metric $\text{LFR}$, defined on $[0,1]$ quantifies the energy ratio of the low-frequency component over the whole spectrum. A higher LFR suggests stronger low-frequency dominance. To avoid high-dimensional Fourier transforms, we utilize Parseval’s identity and the convolution theorem to rewrite LFR as:
\begin{align}
    \text{LFR} (k_0) = \frac{  \displaystyle \int_{\mathbb{R}^d} 
 \left\vert \boldsymbol{y} ( \boldsymbol{x} ) * \mathcal{F}^{-1}[ \mathbbm{1}_{|\boldsymbol{k}|\leq k_0}] \right\vert ^2\mathrm{d} \boldsymbol{x} } { \displaystyle \int_{\mathbb{R}^d} 
 \left\vert  \boldsymbol{y} ( \boldsymbol{x} )  \right\vert ^2\mathrm{d} \boldsymbol{x} } ,
\end{align}
where $\mathcal{F}^{-1}$ represents the inverse Fourier transform and $*$ denotes convolution operation. The inverse transform of the indicator function remains challenging to compute directly. Therefore, we resort to approximating it with a Gaussian function $\hat{G}^{\delta}  (\boldsymbol{k})$ in the frequency domain with variance $1/\delta$. Its inverse Fourier transform in the spatial domain, ${G}^{\delta}(\boldsymbol{x})$, remains a Gaussian with variance $\delta$, possessing a simple and tractable form. Additionally, due to its exponential decay property with $\boldsymbol{k}$, the Gaussian function serves as an effective approximation to $\mathbbm{1}_{|\boldsymbol{k}|\leq k_0}$. Note that the variance $1/\delta$ of $\hat{G}^{\delta}(\boldsymbol{k})$ can be interpreted as $k_0$ in the frequency domain.

We consider the discrete LFR on the dataset $\mathcal{S} := \{ (\boldsymbol{x}_i, \boldsymbol{y}_i) \}_{i=1}^{n} $. The low-frequency component can be expressed as:
\begin{align}
\boldsymbol{y}_i^{\text{low}, \delta\left(k_0\right)} \triangleq\left(\boldsymbol{y} * G^{\delta\left(k_0\right)}\right)_i = \frac{1}{C_i} \sum_{j=1}^{n} \boldsymbol{y}_j G^\delta\left(\boldsymbol{x}_i-\boldsymbol{x}_j\right), 
\end{align}
where \( C_i = \sum_{j=0}^{n-1} G^\delta\left(\boldsymbol{x}_i-\boldsymbol{x}_j\right) \) is a normalization factor, and
\begin{align}
G^\delta\left(\boldsymbol{x}_i-\boldsymbol{x}_j\right) = \exp \left(-\left|\boldsymbol{x}_i-\boldsymbol{x}_j\right |^2 / 2 \delta\right).
\end{align}
Hence, the discrete LFR can be denoted as:
\begin{align}
\operatorname{LFR}\left(k_0\right)=\frac{\sum_i\left|\boldsymbol{y}_i^{\text {low }, \delta\left(k_0\right)}\right|^2}{\sum_i\left|\boldsymbol{y}_i\right|^2} .
\end{align}
$\text{LFR}$ represents the energy ratio within the ball of radius $k_0$ in the frequency domain. Analogous to the probability density function, we further define the Ratio Density Function (RDF) as:
\begin{align}
    \text{RDF}(k_0) = \frac{\partial \text{LFR}(k_0)}{\partial k_0}.
\end{align}
RDF describes the distribution of the low-frequency energy ratio over the spherical surface of radius $k_0$, thereby facilitating indirect visualization and analysis of the spectral properties of high-dimensional functions.

\subsubsection{GBCT scales high frequency to low frequency}

\begin{figure}[htbp]
    \centering
    \includegraphics[width=0.8\linewidth]{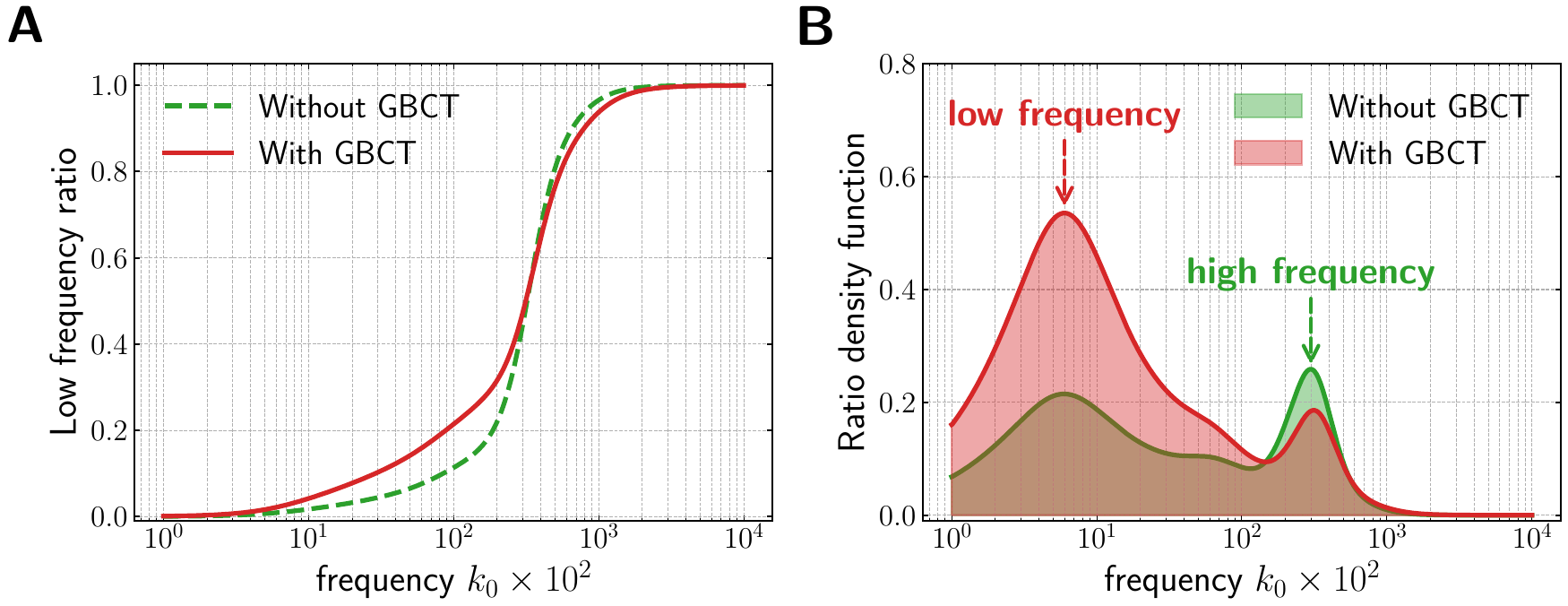}
    \caption{ Frequency analysis for GBCT. (A) Low Frequency Ratio (LFR) of the training data; (B) Ratio Density Function (RDF) of the training data. In the practical calculations, we take $k_0 \approx 1/\delta$ .}
    \label{fig:freq_analysis}
\end{figure}

To investigate the impact of GBCT on the frequency distribution within the training data, we separately computed the LFR and RDF of $ \boldsymbol{u}_B $ and $ G(\boldsymbol{u}_B) $ (defined in Eq.~\ref{eq:uB}), which correspond to the labels of \bctnet{} and \gbctnet{}, respectively. 
In our experiments, we set $ k_0 \approx 1/\delta $. 
The results are illustrated in Fig.~\ref{fig:freq_analysis}. 
After applying GBCT scaling, as shown in Fig.~\ref{fig:freq_analysis}(A), the LFR curve becomes noticeably smoother, indicating a more uniform frequency distribution. Furthermore, Fig.~\ref{fig:freq_analysis}(B) reveals that the RDF of the objective function shifts significantly towards lower frequencies following GBCT scaling, aligning well with the neural network’s inherent preference for low-frequency components.

The low-frequency preference of neural networks has inspired a series of studies. Early works such as MscaleDNN \cite{liu2020CCP}, as well as more recent methods like VLT-PINN \cite{wu2024JCP}, KH-PINN \cite{wu2025PF}, and VS-PINN \cite{ko2025JCP}, have significantly enhanced the performance of PINNs via variable scalings.   
However, when addressing reaction dynamical systems with typical multi-scale nature, data transformation strategies remain under-explored.
Previous deep learning-based algorithms have recognized the necessity of transforming reactant concentrations (e.g., logarithmic transformations and BCT), but they often overlook the subtle small-scale or high-frequency features in target data that are difficult to detect. 
The design of GBCT aims to fill this gap.
As a novel data scaling method, GBCT was initially motivated by the intrinsic properties of elementary reactions in reaction kinetics. 
We emphasize that GBCT bridges neural networks modeling and multiscale dynamical data through physical intuition.

\section{Conclusion}
\label{sec:conclusion}

In this work, we introduce Generalized Box Cox Transformation (GBCT), a novel data scaling method designed for multiscale dynamical systems. We incorporate GBCT as a simple, plug-and-play output-scaling layer within our existing data-driven framework \cite{yao2025CPC}, resulting in two neural networks: \gbctnet{} and the baseline \bctnet{} (without GBCT). 

Comprehensive evaluations across six test scenarios demonstrate the accuracy and superior stability of \gbctnet{}, encompassing 21-species chemical and 13-isotope nuclear reaction kinetics, the typical Robertson problem coupled with diffusion, and practically relevant benchmarks including turbulent chemically reacting diffusion and nuclear reactive flow simulations. In the chemical kinetics model, \gbctnet{} maintains accurate long-term predictions over 50,000 steps with relative errors of $\mathcal{O}(10^{-2})$, which are two orders of magnitude lower than those of \bctnet{}. Similar improvements are observed in the nuclear kinetics and Robertson diffusion problems, where \gbctnet{} effectively suppresses error accumulation and ensures physically consistent evolution. 
In complex chemical and nuclear reactive flow simulations, \gbctnet{} further reduces prediction errors to $20\sim50$\% of the baseline model, accurately capturing flame fronts and propagation speeds in such reaction-diffusion systems. 
These results demonstrate that \gbctnet{} effectively suppresses error accumulation, thereby improving prediction accuracy and ensuring stability during long-term evolution.
Beyond predictive performance, GBCT enhances training efficiency and generalization, achieving comparable accuracy with only about 1/6 of the training epochs required by \bctnet{}.
Frequency analysis indicates that GBCT effectively shifts high-frequency components of the objective function toward lower frequencies, aligning with the neural network’s natural low frequency bias and thereby facilitating model training.

Although many studies tackle multiscale dynamical systems from the perspective of physics knowledge embedding, network architecture or loss function design, we aim to introduce a new data-level method to boost prediction capability through a simple yet effective output scaling.
The innovation of GBCT lies in bridging prior dynamic analysis, which pinpoints the origin of multiscale behaviors, with a targeted scaling transformation that proves highly effective in practice.
Moreover, GBCT is model-agnostic and can be potentially integrated as a plug-in into popular frameworks such as PINNs and operator learning methods to mitigate the multiscale effects arising from sign-changing data across diverse multiscale problems.


\section*{Declaration of Generative AI and AI-assisted technologies in the writing process}

The authors used GPT-4o to polish the language of this manuscript and then thoroughly reviewed and revised the content to ensure its accuracy and completeness. They take full responsibility for the final version.

\section*{Acknowledgment}

This work is sponsored by the National Natural Science Foundation of China Grants No. 92470127(T. Z.), 92270203(T. Z.), the HPC of School of Mathematical Sciences and the Student Innovation Center, and the Siyuan-1 cluster supported by the Center for High Performance Computing at Shanghai Jiao Tong University.

\section*{Compliance with Ethical Standards}
Conflict of Interest: The authors declare that they have no conflict of interest.

\bibliographystyle{elsarticle-num} 
\bibliography{reference} 

\begin{thebibliography}{10}
\expandafter\ifx\csname url\endcsname\relax
  \def\url#1{\texttt{#1}}\fi
\expandafter\ifx\csname urlprefix\endcsname\relax\def\urlprefix{ URL }\fi
\expandafter\ifx\csname href\endcsname\relax
  \def\href#1#2{#2} \def\path#1{#1}\fi

\bibitem{zingale2024}
M.~Zingale, K.~Bhargava, R.~Brady, Z.~Chen, S.~Guichandut, E.~T. Johnson, M.~Katz, A.~S. Clark, The challenges of modeling astrophysical reacting flows (Nov. 2024).
\newblock \href {http://arxiv.org/abs/2411.12491} {\path{arXiv:2411.12491}}, \href {https://doi.org/10.48550/arXiv.2411.12491} {\path{doi:10.48550/arXiv.2411.12491}}.

\bibitem{christo1996CF}
F.~Christo, A.~Masri, E.~Nebot, Artificial neural network implementation of chemistry with pdf simulation of {{H2}}/{{CO2}} flames, Combustion and Flame 106~(4) (1996) 406--427.
\newblock \href {https://doi.org/10.1016/0010-2180(95)00250-2} {\path{doi:10.1016/0010-2180(95)00250-2}}.

\bibitem{christo1996SC}
F.~Christo, A.~Masri, E.~Nebot, S.~Pope, An integrated {{PDF}}/neural network approach for simulating turbulent reacting systems, Symposium (International) on Combustion 26~(1) (1996) 43--48.
\newblock \href {https://doi.org/10.1016/S0082-0784(96)80198-6} {\path{doi:10.1016/S0082-0784(96)80198-6}}.

\bibitem{blasco1998CF}
J.~Blasco, N.~Fueyo, C.~Dopazo, J.~Ballester, Modelling the temporal evolution of a reduced combustion chemical system with an artificial neural network, Combustion and Flame 113~(1) (1998) 38--52.
\newblock \href {https://doi.org/10.1016/S0010-2180(97)00211-3} {\path{doi:10.1016/S0010-2180(97)00211-3}}.

\bibitem{blasco2000CTM}
J.~Blasco, N.~Fueyo, C.~Dopazo, J.-Y. Chen, A self-organizing-map approach to chemistry representation in combustion applications, Combustion Theory and Modelling 4~(1) (2000) 61--76.
\newblock \href {http://arxiv.org/abs/https://doi.org/10.1088/1364-7830/4/1/304} {\path{arXiv:https://doi.org/10.1088/1364-7830/4/1/304}}, \href {https://doi.org/10.1088/1364-7830/4/1/304} {\path{doi:10.1088/1364-7830/4/1/304}}.

\bibitem{kempf2005PCI}
A.~Kempf, F.~Flemming, J.~Janicka, Investigation of lengthscales, scalar dissipation, and flame orientation in a piloted diffusion flame by {{LES}}, Proceedings of the Combustion Institute 30~(1) (2005) 557--565.
\newblock \href {https://doi.org/10.1016/j.proci.2004.08.182} {\path{doi:10.1016/j.proci.2004.08.182}}.

\bibitem{ihme2009PCI}
M.~Ihme, C.~Schmitt, H.~Pitsch, Optimal artificial neural networks and tabulation methods for chemistry representation in {{LES}} of a bluff-body swirl-stabilized flame, Proceedings of the Combustion Institute 32~(1) (2009) 1527--1535.
\newblock \href {https://doi.org/10.1016/j.proci.2008.06.100} {\path{doi:10.1016/j.proci.2008.06.100}}.

\bibitem{raissi2019JCP}
M.~Raissi, P.~Perdikaris, G.~Karniadakis, Physics-informed neural networks: A deep learning framework for solving forward and inverse problems involving nonlinear partial differential equations, Journal of Computational Physics 378 (2019) 686--707.
\newblock \href {https://doi.org/10.1016/j.jcp.2018.10.045} {\path{doi:10.1016/j.jcp.2018.10.045}}.

\bibitem{chen2018ANIP}
R.~T.~Q. Chen, Y.~Rubanova, J.~Bettencourt, D.~K. Duvenaud, Neural ordinary differential equations, in: S.~Bengio, H.~Wallach, H.~Larochelle, K.~Grauman, N.~{Cesa-Bianchi}, R.~Garnett (Eds.), Advances in {{Neural Information Processing Systems}}, Vol.~31, Curran Associates, Inc., 2018.

\bibitem{owoyele2022EA}
O.~Owoyele, P.~Pal, {{ChemNODE}}: A neural ordinary differential equations framework for efficient chemical kinetic solvers, Energy and AI 7 (2022) 100118.
\newblock \href {https://doi.org/10.1016/j.egyai.2021.100118} {\path{doi:10.1016/j.egyai.2021.100118}}.

\bibitem{ji2021JPCA}
W.~Ji, W.~Qiu, Z.~Shi, S.~Pan, S.~Deng, Stiff-{{PINN}}: Physics-informed neural network for stiff chemical kinetics, Journal of Physical Chemistry A 125~(36) (2021) 8098--8106.
\newblock \href {https://doi.org/10.1021/acs.jpca.1c05102} {\path{doi:10.1021/acs.jpca.1c05102}}.

\bibitem{weng2022JPCA}
Y.~Weng, D.~Zhou, Multiscale physics-informed neural networks for stiff chemical kinetics, Journal of Physical Chemistry A 126~(45) (2022) 8534--8543.
\newblock \href {https://doi.org/10.1021/acs.jpca.2c06513} {\path{doi:10.1021/acs.jpca.2c06513}}.

\bibitem{deflorio2022CIJN}
M.~De~Florio, E.~Schiassi, R.~Furfaro, Physics-informed neural networks and functional interpolation for stiff chemical kinetics, Chaos: an Interdisciplinary Journal of Nonlinear Science 32~(6) (2022) 63107.
\newblock \href {https://doi.org/10.1063/5.0086649} {\path{doi:10.1063/5.0086649}}.

\bibitem{lu2021NMI}
L.~Lu, P.~Jin, G.~Pang, Z.~Zhang, G.~E. Karniadakis, Learning nonlinear operators via {{DeepONet}} based on the universal approximation theorem of operators, Nature Machine Intelligence 3~(3) (2021) 218--229.
\newblock \href {https://doi.org/10.1038/s42256-021-00302-5} {\path{doi:10.1038/s42256-021-00302-5}}.

\bibitem{li2020APA}
Z.~Li, N.~Kovachki, K.~Azizzadenesheli, B.~Liu, K.~Bhattacharya, A.~Stuart, A.~Anandkumar, Fourier neural operator for parametric partial differential equations, Arxiv Preprint Arxiv:2010.08895 (2020).
\newblock \href {http://arxiv.org/abs/2010.08895} {\path{arXiv:2010.08895}}.

\bibitem{goswami2024CMAM}
S.~Goswami, A.~D. Jagtap, H.~Babaee, B.~T. Susi, G.~E. Karniadakis, Learning stiff chemical kinetics using extended deep neural operators, Computer Methods in Applied Mechanics and Engineering 419 (2024) 116674.
\newblock \href {https://doi.org/10.1016/j.cma.2023.116674} {\path{doi:10.1016/j.cma.2023.116674}}.

\bibitem{weng2025CF}
Y.~Weng, H.~Li, H.~Zhang, Z.~X. Chen, D.~Zhou, Extended fourier neural operators to learn stiff chemical kinetics under unseen conditions, Combustion and Flame 272 (2025) 113847.
\newblock \href {https://doi.org/10.1016/j.combustflame.2024.113847} {\path{doi:10.1016/j.combustflame.2024.113847}}.

\bibitem{wan2020CF}
K.~Wan, C.~Barnaud, L.~Vervisch, P.~Domingo, Chemistry reduction using machine learning trained from non-premixed micro-mixing modeling: Application to {{DNS}} of a syngas turbulent oxy-flame with side-wall effects, Combustion and Flame 220 (2020) 119--129.
\newblock \href {https://doi.org/10.1016/j.combustflame.2020.06.008} {\path{doi:10.1016/j.combustflame.2020.06.008}}.

\bibitem{ding2021CF}
T.~Ding, T.~Readshaw, S.~Rigopoulos, W.~Jones, Machine learning tabulation of thermochemistry in turbulent combustion: An approach based on hybrid flamelet/random data and multiple multilayer perceptrons, Combustion and Flame 231 (2021) 111493.
\newblock \href {https://doi.org/10.1016/j.combustflame.2021.111493} {\path{doi:10.1016/j.combustflame.2021.111493}}.

\bibitem{readshaw2023CF}
T.~Readshaw, L.~L.~C. Franke, W.~P. Jones, S.~Rigopoulos, Simulation of turbulent premixed flames with machine learning - tabulated thermochemistry, Combustion and Flame 258 (2023) 113058.
\newblock \href {https://doi.org/10.1016/j.combustflame.2023.113058} {\path{doi:10.1016/j.combustflame.2023.113058}}.

\bibitem{zhang2022CF}
T.~Zhang, Y.~Yi, Y.~Xu, Z.~X. Chen, Y.~Zhang, W.~E, Z.-Q.~J. Xu, A multi-scale sampling method for accurate and robust deep neural network to predict combustion chemical kinetics, Combustion and Flame 245 (2022) 112319.
\newblock \href {https://doi.org/10.1016/j.combustflame.2022.112319} {\path{doi:10.1016/j.combustflame.2022.112319}}.

\bibitem{yao2025CPC}
J.~Yao, Y.~Yi, L.~Hang, W.~E, W.~Wang, Y.~Zhang, T.~Zhang, Z.-Q.~J. Xu, Solving multiscale dynamical systems by deep learning, Computer Physics Communications 316 (2025) 109802.
\newblock \href {https://doi.org/10.1016/j.cpc.2025.109802} {\path{doi:10.1016/j.cpc.2025.109802}}.

\bibitem{zhang2025AJ}
X.~Zhang, Y.~Yi, L.~Wang, Z.-Q.~J. Xu, T.~Zhang, Y.~Zhou, Deep neural networks for modeling astrophysical nuclear reacting flows, Astrophysical Journal 990~(2) (2025) 105.
\newblock \href {https://doi.org/10.3847/1538-4357/adf331} {\path{doi:10.3847/1538-4357/adf331}}.

\bibitem{wang2025CF}
T.~Wang, Y.~Yi, J.~Yao, Z.-Q.~J. Xu, T.~Zhang, Z.~Chen, Enforcing physical conservation in neural network surrogate models for complex chemical kinetics, Combustion and Flame 275 (2025) 114105.
\newblock \href {https://doi.org/10.1016/j.combustflame.2025.114105} {\path{doi:10.1016/j.combustflame.2025.114105}}.

\bibitem{Xu2020CiCP}
Z.-Q. J.~X. Zhi-Qin John~Xu, Y.~Z. Yaoyu~Zhang, T.~L. Tao~Luo, Y.~X. Yanyang~Xiao, Z.~M. Zheng~Ma, \href{https://global-sci.com/article/79739/frequency-principle-fourier-analysis-sheds-light-on-deep-neural-networks}{Frequency {{Principle}}: {{Fourier Analysis Sheds Light}} on {{Deep Neural Networks}}} 28~(5)  1746--1767.
\newblock \href {https://doi.org/10.4208/cicp.OA-2020-0085} {\path{doi:10.4208/cicp.OA-2020-0085}}.

\bibitem{xu2025CAMC}
Z.-Q.~J. Xu, Y.~Zhang, T.~Luo, Overview {{Frequency Principle}}/{{Spectral Bias}} in {{Deep Learning}}, Communications on Applied Mathematics and Computation 7~(3) (2025) 827--864.
\newblock \href {https://doi.org/10.1007/s42967-024-00398-7} {\path{doi:10.1007/s42967-024-00398-7}}.

\bibitem{xu2019NIP}
Z.-Q.~J. Xu, Y.~Zhang, Y.~Xiao, \href{http://link.springer.com/10.1007/978-3-030-36708-4_22}{Training behavior of deep neural network in frequency domain}, in: T.~Gedeon, K.~W. Wong, M.~Lee (Eds.), Neural {{Information Processing}}, Vol. 11953, Springer International Publishing, pp. 264--274.
\newblock \href {https://doi.org/10.1007/978-3-030-36708-4_22} {\path{doi:10.1007/978-3-030-36708-4_22}}.

\bibitem{rahaman19a}
N.~Rahaman, A.~Baratin, D.~Arpit, F.~Draxler, M.~Lin, F.~Hamprecht, Y.~Bengio, A.~Courville, \href{https://proceedings.mlr.press/v97/rahaman19a.html}{On the spectral bias of neural networks}, in: K.~Chaudhuri, R.~Salakhutdinov (Eds.), Proceedings of the 36th International Conference on Machine Learning, Vol.~97 of Proceedings of Machine Learning Research, PMLR, 2019, pp. 5301--5310.

\bibitem{kim2021CIJN}
S.~Kim, W.~Ji, S.~Deng, Y.~Ma, C.~Rackauckas, Stiff neural ordinary differential equations, Chaos: an Interdisciplinary Journal of Nonlinear Science 31~(9) (2021) 93122.
\newblock \href {https://doi.org/10.1063/5.0060697} {\path{doi:10.1063/5.0060697}}.

\bibitem{tz2021}
T.~Zhang, Y.~Zhang, W.~E, Y.~Ju, \href{https://arc.aiaa.org/doi/abs/10.2514/6.2021-1139}{DLODE: a deep learning-based ODE solver for chemistry kinetics}, 2021.
\newblock \href {http://arxiv.org/abs/https://arc.aiaa.org/doi/pdf/10.2514/6.2021-1139} {\path{arXiv:https://arc.aiaa.org/doi/pdf/10.2514/6.2021-1139}}, \href {https://doi.org/10.2514/6.2021-1139} {\path{doi:10.2514/6.2021-1139}}.

\bibitem{lian2023AA}
C.~Lian, T.~Tang, H.~Wang, J.~Yu, M.~Sun, D.~Xiong, Y.~Yang, Study on the application of artificial neural network-based flamelet/progress variable model in supersonic combustion, AIP Advances 13~(11) (2023) 115319.
\newblock \href {https://doi.org/10.1063/5.0171442} {\path{doi:10.1063/5.0171442}}.

\bibitem{saito2023AECS}
M.~Saito, J.~Xing, J.~Nagao, R.~Kurose, Data-driven simulation of ammonia combustion using neural ordinary differential equations ({{NODE}}), Applications in Energy and Combustion Science 16 (2023) 100196.
\newblock \href {https://doi.org/10.1016/j.jaecs.2023.100196} {\path{doi:10.1016/j.jaecs.2023.100196}}.

\bibitem{lei2024ANE}
K.~Lei, H.~Wu, Z.~Liu, Y.~Cao, G.~Liu, X.~Li, Q.~He, L.~Cao, {{SN-MscaleDNN}}: A coupling approach for rapid shielding-scheme evaluation of micro gas-cooled reactor in the large design-parameter space, Annals of Nuclear Energy 196 (2024) 110241.
\newblock \href {https://doi.org/10.1016/j.anucene.2023.110241} {\path{doi:10.1016/j.anucene.2023.110241}}.

\bibitem{wu2024EA}
S.~Wu, H.~Wang, K.~H. Luo, A robust autoregressive long-term spatiotemporal forecasting framework for surrogate-based turbulent combustion modeling via deep learning, Energy and AI 15 (2024) 100333.
\newblock \href {https://doi.org/10.1016/j.egyai.2023.100333} {\path{doi:10.1016/j.egyai.2023.100333}}.

\bibitem{zhang2024PF}
M.~Zhang, R.~Mao, H.~Li, Z.~An, Z.~X. Chen, Graphics processing unit/artificial neural network-accelerated large-eddy simulation of swirling premixed flames, Physics of Fluids 36~(5) (2024) 55147.
\newblock \href {https://doi.org/10.1063/5.0202321} {\path{doi:10.1063/5.0202321}}.

\bibitem{li2025}
H.~Li, R.~Yang, Y.~Xu, M.~Zhang, R.~Mao, Z.~X. Chen, Comprehensive deep learning for combustion chemistry integration: Multi-fuel generalization and a posteriori validation in reacting flow, Physics of Fluids 37~(1) (2025) 015162.
\newblock \href {https://doi.org/10.1063/5.0248582} {\path{doi:10.1063/5.0248582}}.

\bibitem{rao2025CF}
S.~Rao, W.~Zhou, W.~Han, Y.~Tang, X.~Xu, An adaptive implicit time-integration scheme for stiff chemistry based on jacobian tabulation method, Combustion and Flame 274 (2025) 113997.
\newblock \href {https://doi.org/10.1016/j.combustflame.2025.113997} {\path{doi:10.1016/j.combustflame.2025.113997}}.

\bibitem{robert1966}
H.~Robertson, The solution of a set of reaction rate equations, in: J.~Walsh (Ed.), Numerical Analysis: An Introduction, Vol. 178182, Academic Press, London, 1966, pp. 178--182.

\bibitem{smith2023A}
A.~I. Smith, E.~T. Johnson, Z.~Chen, K.~Eiden, D.~E. Willcox, B.~Boyd, L.~Cao, C.~J. DeGrendele, M.~Zingale, Pynucastro: {{A Python Library}} for {{Nuclear Astrophysics}}, The Astrophysical Journal 947~(2) (2023) 65.
\newblock \href {https://doi.org/10.3847/1538-4357/acbaff} {\path{doi:10.3847/1538-4357/acbaff}}.

\bibitem{Strang1968}
G.~Strang, \href{https://doi.org/10.1137/0705041}{On the construction and comparison of difference schemes}, SIAM Journal on Numerical Analysis 5~(3) (1968) 506--517.
\newblock \href {http://arxiv.org/abs/https://doi.org/10.1137/0705041} {\path{arXiv:https://doi.org/10.1137/0705041}}, \href {https://doi.org/10.1137/0705041} {\path{doi:10.1137/0705041}}.

\bibitem{wang2021}
S.~Wang, Y.~Teng, P.~Perdikaris, Understanding and mitigating gradient flow pathologies in physics-informed neural networks, SIAM Journal on Scientific Computing 43~(5) (2021) A3055--A3081.

\bibitem{Box1964}
G.~E.~P. Box, D.~R. Cox, \href{https://rss.onlinelibrary.wiley.com/doi/abs/10.1111/j.2517-6161.1964.tb00553.x}{An analysis of transformations}, Journal of the Royal Statistical Society: Series B (Methodological) 26~(2) (1964) 211--243.
\newblock \href {http://arxiv.org/abs/https://rss.onlinelibrary.wiley.com/doi/pdf/10.1111/j.2517-6161.1964.tb00553.x} {\path{arXiv:https://rss.onlinelibrary.wiley.com/doi/pdf/10.1111/j.2517-6161.1964.tb00553.x}}, \href {https://doi.org/https://doi.org/10.1111/j.2517-6161.1964.tb00553.x} {\path{doi:https://doi.org/10.1111/j.2517-6161.1964.tb00553.x}}.

\bibitem{rivas_impact_2022}
F.~Rivas, J.~A. Harris, W.~R. Hix, O.~E.~B. Messer, \href{https://iopscience.iop.org/article/10.3847/1538-4357/ac8b06}{The {Impact} of {Resolution} on {Double}-detonation {Models} for {Type} {Ia} {Supernovae}}, The Astrophysical Journal 937~(1) (2022) 2, number: 1.
\newblock \href {https://doi.org/10.3847/1538-4357/ac8b06} {\path{doi:10.3847/1538-4357/ac8b06}}.

\bibitem{saad2017}
T.~Saad, D.~Cline, R.~Stoll, J.~C. Sutherland, \href{https://doi.org/10.2514/1.J055230}{Scalable tools for generating synthetic isotropic turbulence with arbitrary spectra}, AIAA Journal 55~(1) (2017) 327--331.
\newblock \href {https://doi.org/10.2514/1.J055230} {\path{doi:10.2514/1.J055230}}.

\bibitem{xu2021AAAI}
Z.~J. Xu, H.~Zhou, Deep frequency principle towards understanding why deeper learning is faster, Proceedings of the AAAI Conference on Artificial Intelligence 35~(12) (2021) 10541--10550.
\newblock \href {https://doi.org/10.1609/aaai.v35i12.17261} {\path{doi:10.1609/aaai.v35i12.17261}}.

\bibitem{liu2020CCP}
W.~Liu, ZiqiCai, Z.-Q. John~Xu, Multi-scale deep neural network ({{MscaleDNN}}) for solving poisson-boltzmann equation in complex domains, Communications in Computational Physics 28~(5) (2020) 1970--2001.
\newblock \href {https://doi.org/10.4208/cicp.OA-2020-0179} {\path{doi:10.4208/cicp.OA-2020-0179}}.

\bibitem{wu2024JCP}
J.~Wu, Y.~Wu, G.~Zhang, Y.~Zhang, Variable linear transformation improved physics-informed neural networks to solve thin-layer flow problems, Journal of Computational Physics 500 (2024) 112761.
\newblock \href {https://doi.org/10.1016/j.jcp.2024.112761} {\path{doi:10.1016/j.jcp.2024.112761}}.

\bibitem{wu2025PF}
J.~Wu, Y.~Wu, X.~Li, G.~Zhang, Physics-informed neural networks for kelvin--helmholtz instability with spatiotemporal and magnitude multiscale, Physics of Fluids 37~(3) (2025) 34118.
\newblock \href {https://doi.org/10.1063/5.0251167} {\path{doi:10.1063/5.0251167}}.

\bibitem{ko2025JCP}
S.~Ko, S.~Park, {{VS-PINN}}: A fast and efficient training of physics-informed neural networks using variable-scaling methods for solving {{PDEs}} with stiff behavior, Journal of Computational Physics 529 (2025) 113860.
\newblock \href {https://doi.org/10.1016/j.jcp.2025.113860} {\path{doi:10.1016/j.jcp.2025.113860}}.

\bibitem{zirwes2020FTCa}
T.~Zirwes, F.~Zhang, P.~Habisreuther, M.~Hansinger, H.~Bockhorn, M.~Pfitzner, D.~Trimis, Quasi-{{DNS}} dataset of a piloted flame with inhomogeneous inlet conditions, Flow, Turbulence and Combustion 104~(4) (2020) 997--1027.
\newblock \href {https://doi.org/10.1007/s10494-019-00081-5} {\path{doi:10.1007/s10494-019-00081-5}}.

\bibitem{zhang2015HPCS}
F.~Zhang, H.~Bonart, T.~Zirwes, P.~Habisreuther, H.~Bockhorn, N.~Zarzalis, Direct numerical simulation of chemically reacting flows with the public domain code {{OpenFOAM}}, in: W.~E. Nagel, D.~H. Kr{\"o}ner, M.~M. Resch (Eds.), High {{Performance Computing}} in {{Science}} and {{Engineering}} `14, Springer International Publishing, Cham, 2015, pp. 221--236.

\end{thebibliography}

\newpage

\begin{appendices}

\part{}
\vspace{-10em}

\setcounter{parttocdepth}{3}
\part{Appendix} 
\parttoc 

\section{Additional results} 

\subsection{Robertson diffusion problem}
\label{sec:robert_appd}

\begin{figure}[htbp]
	\centering
	\includegraphics[width=1.0\linewidth]{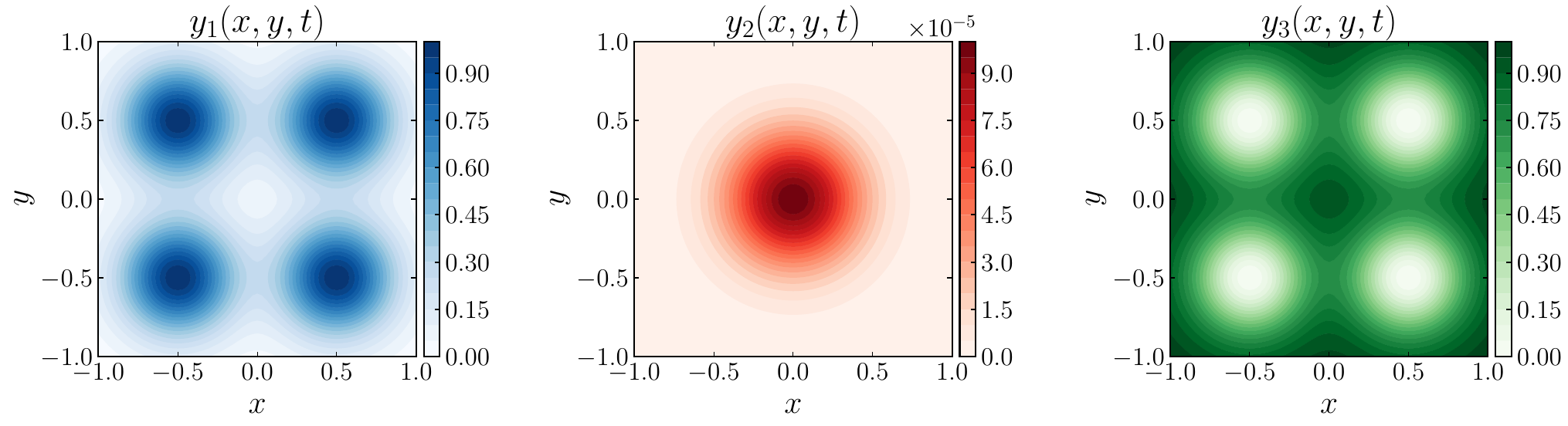}
	\caption{Robertson diffusion initial setup. }
	\label{fig:robert_2d_init}
\end{figure}

The computational domain is configured as $[-1,1]\times[-1,1]$, discretized using a grid with $N_x=N_y=64$ points. Initially, $y_1$ and $y_2$ are generated according to Gaussian functions, while $y_3$ is determined by the constraint condition $\sum_{i=1}^{3} y_i=1$:
\begin{align}
y_i(x,y, t=0) = A_i\cdot \exp\left( - \frac{(x-c_x)^2+(y-c_y)^2}{2\sigma_i^2} \right), \quad i=1,2
\end{align}
where $A_1=1.0$, $\sigma_1=0.25$, $A_2=1\times10^{-4}$, and $\sigma_2 = 0.3$. 
The coordinates $(c_x, c_y)$ represent the center of the Gaussian functions. 
For $y_1$, the center points are located at $(-0.5, 0.5)$, $(-0.5, -0.5)$, $(0.5, -0.5)$, and $(0.5, 0.5)$, respectively. For $y_2$, the center point i s positioned at $(0,0)$. Then we define $y_3=1-y_1-y_2$ and enforce the physical constraints, such as non-negativity and summation to unity, via clipping and normalization.

\subsection{Governing equations of DRM19 reaction-diffusion systems}
\label{sec:ebi_equation_appd}
In this study, the chemical kinetics model used here is DRM19, a reduced reactions subset of full GRI-Mech 1.2, including 21 species and 84 reactions~\footnote{A. Kazakov, M. Frenklach, Reduced Reaction Sets based on GRI-Mech 1.2, \url{http://www.me.berkeley.edu/drm/}.}. We utilize the \textit{EBI-DNS} code \cite{zirwes2020FTCa, zhang2015HPCS} to conduct the numerical simulation of compressible reactive flow involving $N$ chemical species using Direct Numerical Simulation (DNS). The conservation equations of mass, momentum, species and energy can be expressed as follows: 
\begin{align}
\frac{\partial \rho}{\partial t} + \nabla \cdot (\rho \boldsymbol{v}) &= 0~, \\
\frac{\partial (\rho \boldsymbol{v}) }{\partial t} + \nabla \cdot (\rho \boldsymbol{v} \boldsymbol{v}) &= -\nabla p + \nabla \cdot \boldsymbol{\tau} + \rho \mathbf{g}~, \\
\frac{\partial (\rho Y_k)}{\partial t} + \nabla \cdot (\rho Y_k \boldsymbol{v}) &= -\nabla \cdot \mathbf{j}_k + \dot{\omega}_k, \quad k = 1, 2, \ldots, N~, \\
\frac{\partial \left(\rho \left(h_s+\frac{1}{2}\boldsymbol{v}\cdot\boldsymbol{v}\right)\right) }{\partial t}  &= - \nabla \cdot \left(\rho\boldsymbol{v} (h_s+\frac{1}{2}\boldsymbol{v}\cdot\boldsymbol{v} )\right) -\nabla \cdot \mathbf{q} + \frac{\partial p}{\partial t} + \dot{q}_r~, \\
p &= \rho RT~,
\end{align}
where $\rho$ represents the density, $\boldsymbol{v}$ is the velocity, $p$ and $T$ correspond to the static pressure and temperature, respectively, and $h_s$ denotes the sensible enthalpy. 
The variables $Y_k$ and $\dot{\omega}_k$ represent the mass fraction and chemical reaction rate for species $k$, respectively. 
$R$ stands for the specific gas constant. 
The gravitational body force $\rho \mathbf{g}$ serves as an external force acting on the fluid element. 
In this formulation, heat sources arising from viscous dissipation and radiative effects are assumed negligible.

The mixture-averaged transport model is employed to compute the species diffusion flux $\mathbf{j}_k$, the viscous stress tensor $\boldsymbol{\tau}$ for Newtonian fluids, and the conductive heat flux $\mathbf{q}$:
\begin{align}
\boldsymbol{\tau} &= -\mu\left[\nabla \boldsymbol{v} + (\nabla \boldsymbol{v})^T\right] + \frac{2}{3}\mu (\nabla \cdot \boldsymbol{v})\mathbf{I}, \\
\mathbf{j}_k &= -\rho D_k \nabla Y_k, \\
\mathbf{q} &= -\lambda \nabla T + \sum_{k=1}^{N} \mathbf{j}_k h_k,
\end{align}
where $D_k$ denotes the mixture-averaged diffusion coefficient relating the $k$-th species to the background mixture. The parameters $\mu$ and $\lambda$ represent the mixture-averaged dynamic viscosity and thermal conductivity, respectively. $h_k$ signifies the specific enthalpy of species $k$. The chemical reaction rate $\dot{\omega}_k$ appearing in equation (3) is determined through chemical kinetics rate laws combined with the extended Arrhenius formulation. The energy release from chemical reactions contributes a source term to the sensible enthalpy equation (4):
\begin{align}
\dot{q}_r = -\sum_{k=1}^{N} \dot{\omega}_k h_k^0 
\end{align}
where $h_k^0$ represents the standard formation enthalpy of species $k$.


\subsection{One-dimensional steady laminar flame}
\label{sec:case3_1d_chem_appd}

\begin{figure}[htbp]
    \centering
    \includegraphics[width=1.0\linewidth]{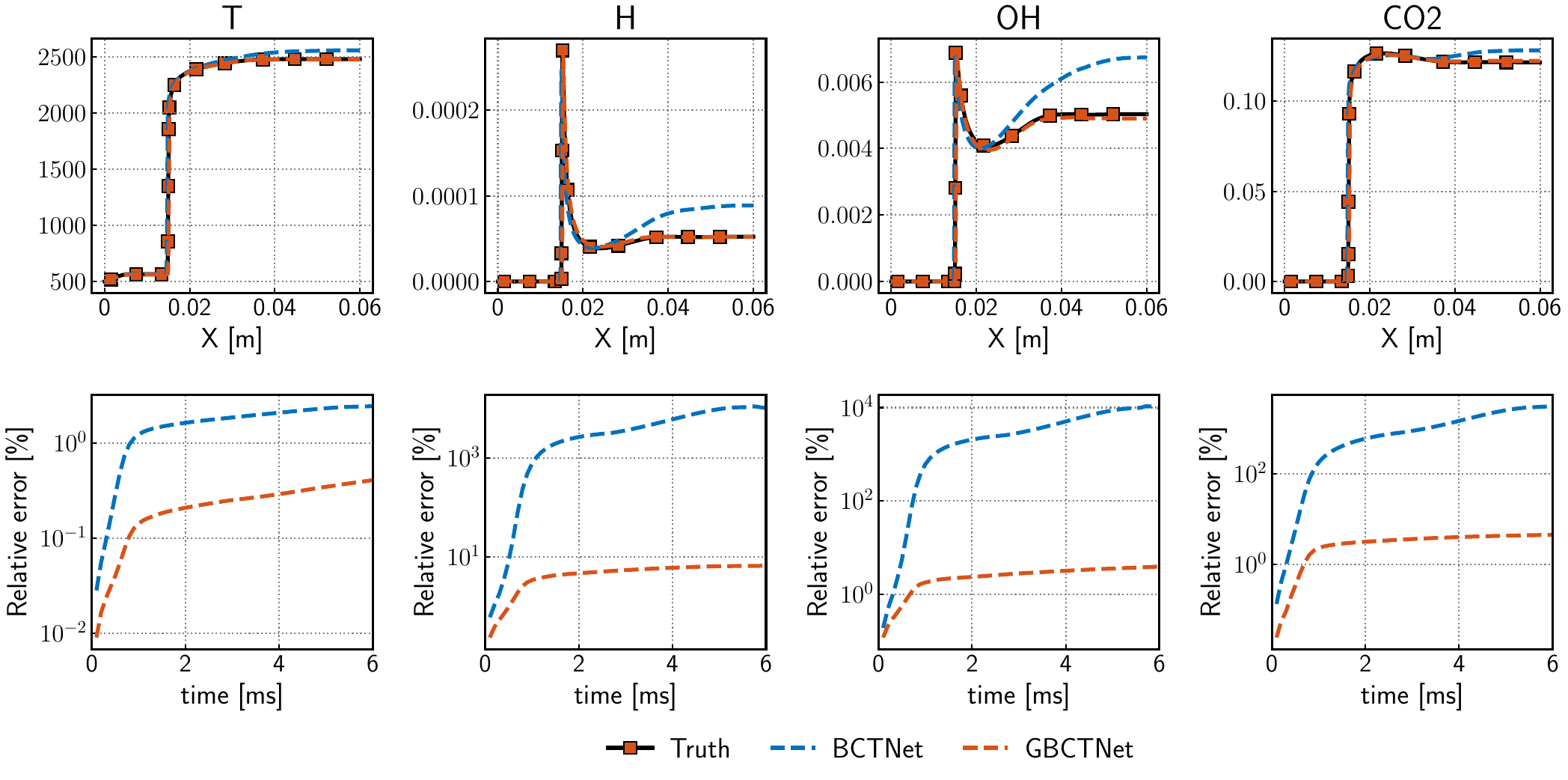}
    \caption{Temperature and mass fraction of radicals \ce{H}, \ce{OH}, \ce{CO2} distribution (upper row) at time = 6 ms and domain-average relative error propagation over time (lower row). The initial temperature $T_0 = 500$K, pressure $p_0=1$ atm and equivalence ratio $\phi_0 =1.0$. }
    \label{fig:drm19_1d}
\end{figure}

We begin by examining a one-dimensional case of a laminar steady-state flame governed by chemical kinetics. The setup involves a duct aligned along the x-axis, characterized by inlet and outlet boundaries. The computational domain spans a length of 6 cm, discretized into 4030 grid cells to finely resolve the flame front structure. The initial conditions are set as follows: temperature $T_0 = 500$ K, pressure $p_0 = 1$ atm, and equivalence ratio $\phi_0 = 1.0$. Using these initial settings, we first compute the steady-state flame solutions. Subsequently, the temporal evolution over 6 ms is simulated, employing DNN models as the ODE solver to assess their accuracy and stability in predicting coupled flow-reaction dynamics.

As depicted in the upper row of Fig.~\ref{fig:drm19_1d}, the spatial distributions of temperature $T$ and the mass fractions of radicals \ce{H}, \ce{OH}, and \ce{CO2} at time = 6 ms are compared. It is observed that \gbctnet{} accurately captures the distribution of physical quantities, whereas \bctnet{} exhibits deviations in the high-temperature region. Furthermore, as shown in the lower row of Fig.~\ref{fig:drm19_1d}, the propagation of average relative error over time is presented. The maximum domain-average relative error in temperature over time predicted by \gbctnet{} is 0.4\%, significantly outperforming the 2.1\% error predicted by \bctnet{}. Additionally, the domain-average relative error in radical mass fractions predicted by \gbctnet{} is an order of magnitude smaller than that predicted by \bctnet{}, highlighting \gbctnet{}'s superior capability in mitigating error accumulation.

\end{appendices}

\end{document}